\documentclass[12pt,a4paper]{amsart}

\usepackage[american]{babel}
\usepackage{bm,amsmath,amssymb}
\usepackage{amsbsy,amsfonts,array}
\usepackage{verbatim,graphicx}
\usepackage{epstopdf}
\usepackage{subfigure}
\graphicspath{{./Figures/},{./}}
\usepackage{color}

\newcommand{\Div}{\text{div}}

\newcommand{\Frac}[2]{\displaystyle\frac{#1}{#2}}
\newcommand{\Sum}[2]{\displaystyle{\sum\limits_{#1}^{#2}}}
\newcommand{\vect}[1]{\bm{#1}}

\newcommand{\unit}[1]{\mathrm{#1}}

\newcommand{\meter}{\unit{m}}
\newcommand{\second}{\unit{s}}
\newcommand{\volt}{\unit{V}}

\newtheorem{remark}{Remark}[section]
\newtheorem{example}{Example}[section]
\newtheorem{proposition}{Proposition}[section]

\title[Electro-Thermo-Chemical 3D Models]
{Electro-Thermo-Chemical Computational Models for 3D Heterogeneous 
Semiconductor Device Simulation}

\address{$^1$ Modeling and Simulations group
         Micron Technology 
         via Olivetti 2, Agrate Brianza, Italy \\
				 $^{2}$ Dipartimento di Matematica \lq\lq F. Brioschi\rq\rq,
         Politecnico di Milano 
				 \\Piazza L. da Vinci 32, 20133 Milano, Italy}
\author{A. Mauri$^{1}$ \and R. Sacco$^{2}$ \and M. Verri$^{2}$}
\email{aureliogianc@micron.com} 
\email{riccardo.sacco@polimi.it}
\email{maurizio.verri@polimi.it}

\begin{document}

\date{\today}

\begin{abstract}
In this article we propose and numerically implement 
a mathematical model for the simulation of three-dimensional
semiconductor devices characterized by an heterogeneous 
material structure. The model consists of a system of 
nonlinearly coupled time-dependent dif\-fu\-sion-re\-ac\-tion 
partial differential equations with convection terms
describing the principal electrical, 
thermal and chemical phenomena that determine the macroscopic electrical
response of the device under the action of externally applied
electrical and thermal forces. The system is supplied with suitable initial,
boundary and interface conditions that account for the interaction
occurring among the various regions of the device with the surrounding
environment. Temporal semi-discretization of the problem is carried 
out with the Backward Euler Method while a fixed-point iteration
of Gummel type is used for system decoupling. Numerical approximation
of the linearized subproblems is carried out using an exponentially 
fitted stabilized Finite Element Method on unstructured tetrahedral grids. 
Several computational experiments are included to validate the
physical accuracy of the proposed computational algorithm in the study
of realistic device structures.
\end{abstract}

\maketitle


\section{Introduction}\label{sec:intro}

The continuous scaling of semiconductor devices 
has pushed contemporary research and most prominent technologies
towards the use of innovative materials where new physical phenomena occur. 
In this context, an important class of applications is represented by 
resistive memories. In the case of Phase Change Memories (PCM) the resistive 
state is determined by a controlled switch of a calchogenide between the 
crystalline and the amorphous phase~\cite{pcm2010}.
Here recent studies have clearly demonstrated the onset of a significant 
mass transport among different components of the calchogenide alloys. 
In other devices, the Resistive Random Access Memories (ReRAM),  
the low and high resistance state~\cite{rram2011} 
is realized by using and controlling 
non-equilibrium thermo-chemical reactions
Moreover, in most of these new applications, 
the active material of the device (where transport, 
diffusion and reaction processes occur) is no longer homogeneous but
often displays a markedly heterogeneous structure, 
as in the case of advanced logic devices.
Finally, during the specific device application, the 
main physical material properties are not constant
but also evolve in time due to the extreme working
conditions (i.e., high electric and/or thermal fields).

A multidisciplinary approach is clearly fundamental
to describe the basic functionality of 
heterogeneous devices in the correct physical framework. 
As a matter of fact, even if the exploration of materials properties
can be effectively understood and theoretically simulated 
with the help of \textquotedblleft ab-initio~\textquotedblright 
calculations~\cite{car1985}, the electrical response
and the time scale of operation of such devices still need to be 
addressed with the advanced mathematical methods 
traditionally employed in electronic, mechanical and thermal simulation. 
The novel challenge introduced by the technological application
considered in the present article is that theoretical 
elements of semiconductor device physics, chemical, thermal 
and mechanical properties, must be included within a unified model 
setting in order to allow self-consistent calculations
that account for the mutual interplay among the various 
phenomena occuring in the same device. This strong requirement
reflects into a similar constraint in the numerical treatment of
the problem because standard simulation suites are no longer usable
but they need to be integrated and in some cases 
completely developed from scratch.

For these reasons, in this article we have developed a 
general mathematical and numerical framework
in which the different physical contributions to the simulation
can be effectively incorporated and mutually coupled to reach
the desidered self-consistency and model accuracy. 

The mathematical model consists of a system of 
nonlinearly coupled time-dependent dif\-fu\-sion-re\-ac\-tion 
partial differential equations (PDEs) with convection terms
describing the principal electrical, thermal and chemical phenomena 
that determine the macroscopic electrical
response of the device under the action of externally applied
electrical and thermal forces (see~\cite{bardfalkner2001,degrootmazur1984}
and~\cite{juengel1997,glitzky1997,glitzky2013,consiglieri2013}). 
The system is supplied with initial,
boundary and interface conditions that account for the interaction
occurring among the various regions of the device with the surrounding
environment. 

The numerical approximation of the problem is conducted in 
two distinct steps. In the first step, temporal semi-discretization 
is carried 
out with the Backward Euler Method using a non-uniform time stepping.
In the second step, a fixed-point iteration of Gummel 
type is adopted for system decoupling~\cite{Jerome1996}. 
This leads to solving a sequence of 
linearized advection-diffusion-reaction equations 
that are numerically treated using an exponentially 
fitted stabilized Finite Element Method 
(FEM)~\cite{Gatti1998,XuZikatanov1999,CdFThesis2006} on 
unstructured tetrahedral partition of the computational domain. 
The FEM is chosen in the present discrete formulation of our model 
because it can properly address 
the complexity of the three-dimensional geometry (3D),
avoiding any requirement of symmetry often used as a simplification 
and offering at the same time the adequate flexibility to implement all the 
mathematical and physical models needed in these emerging applications.

An outline of the article is as follows.
Sect.~\ref{sec:th_el_ch} illustrates the fundamental conservation
laws that express mass and energy balance of a system 
of $M$ charged species in a material medium under the 
combined effect of electrical, thermal and chemical forces.
Sect.~\ref{sec:geometry} is devoted to the 
description of the multi-domain geometrical structure of the 
3D semiconductor device object of the present study while
Sect.~\ref{sec:model} describes how to adapt 
the general thermo-electrochemical theory 
of Sect.~\ref{sec:th_el_ch} to the mathematical modeling of 
the class of devices of Sect.~\ref{sec:geometry}.
The resulting formulation deals with the case of a single negatively 
charged species (electrons, $M=1$) and consists of a nonlinearly 
coupled system of advection-diffusion-reaction PDEs 
that have to be solved in a heterogeneous domain supplied by
a set of initial and boundary conditions.
Sect.~\ref{sec:comp_tech} is, instead, devoted to illustrate the 
three main computational steps which allow to translate
the differential problem of Sect.~\ref{sec:model} into the 
successive solution of linear algebraic systems providing 
the approximate solution of the problem. 
Sect.~\ref{sec:simulations} is devoted to the validation 
of the physical accuracy of the computational model 
through the simulation of 3D device structures under realistic
working conditions. Sect.~\ref{sec:conclusions} 
draws the main conclusions reached in the present article and 
addresses possible future research developments. 
Appendix~\ref{sec:list} contains
a list of all the symbols introduced in the article, specifying
for each symbol the associated physical meaning and units.

\section{Modeling of Thermo-Electrochemical Phenomena}
\label{sec:th_el_ch}

In this section we introduce the fundamental conservation
laws that express mass and energy balance of a charged multi-species
system moving in a material medium under the combined effect of electrical, 
thermal and chemical forces. For a complete treatment of
electrochemical phenomena and of the mathematical foundations
of non-equilibrium thermodynamics, we refer to~\cite{bardfalkner2001}
and to, e.g.,~\cite{degrootmazur1984}. 
For the mathematical analysis of general reaction-diffusion 
thermo-chemically coupled systems, we refer, e.g., 
to~\cite{juengel1997,glitzky1997,glitzky2013,consiglieri2013} 
and to the bibliography cited therein.

Let $M \geq 1$ be the total number of chemicals flowing 
in the medium under the
action of electrical, chemical and thermal forces. 
We denote by $N_i=N_i(\vect{x},t)$, $i=1, \ldots, M$, 
the number density of the $i$-th chemical 
at the spatial position $\vect{x}$ and time $t$, 
and by $z_i$ its ionic valence (equal to zero if 
the species is electrically neutral). 
We set $\mathbf{N}:= \left[ N_1, \ldots, N_M \right]^T$. 
We also introduce the dependent variables $T=T(\vect{x}, t)$ 
and $\vect{E}=\vect{E}(\vect{x}, t)$ representing the 
temperature of the medium and the electric field 
at the spatial position $\vect{x}$ and time $t$,
respectively. 

\subsection{Conservation laws}\label{sec:cons_laws}
The basic form of the mathematical model 
considered in this article is constituted by the following 
coupled system of PDEs in conservation form:
\begin{subequations}\label{eq:basic_system}
\begin{align}
q z_i \Frac{\partial N_i}{\partial t} + 
\Div \vect{j}_i(N_i,T, \vect{E}) = q R_i(\mathbf{N},T, \vect{E})
& \qquad i=1, \ldots, M \label{eq:mass_balance_N_i} \\[2mm]
\Frac{\partial}{\partial t}  (\rho c T) + 
\Div \vect{j}_T(\mathbf{N},T, \vect{E}) = \mathcal{Q}_T(\mathbf{N},T,
\vect{E}) &
\label{eq:energy_balance}\\[2mm]
\Div(\varepsilon \vect{E}) = q\mathcal{D} + \Sum{i=1}{M} q z_i N_i 
& \label{eq:Poisson}
\end{align}

Eqns.~\eqref{eq:mass_balance_N_i} are the continuity
equations for the $M$ chemicals $N_i$, $i=1, \ldots, M$, 
where $\vect{j}_i$ is the current density associated with each chemical $N_i$
and $R_i$ is the corresponding net production rate accounting
for recombination and generation phenomena in the medium.

Eq.~\eqref{eq:energy_balance} is the energy balance equation in
the system, where $\rho$ and $c$ are the mass density and the
specific heat of the medium,
respectively, $\vect{j}_T$ is the energy flux density in the medium
while $\mathcal{Q}_T$ is the net heat production rate. 

Eq.~\eqref{eq:Poisson} is the Poisson equation expressing
Gauss' law in differential form, where $q$ is the electron charge,
$\varepsilon$ is the dielectric permittivity of the medium
and $\mathcal{D}$ is a given function of position 
that accounts for the possible presence of fixed ionized dopant impurities.
Assuming the quasi-static approximation in Maxwell's equations
(see~\cite{Selberherr84}), the electric field can be expressed as
\begin{equation}\label{eq:efield}
\vect{E} = -\nabla \varphi
\end{equation}
\end{subequations}
where $\varphi=\varphi(\vect{x},t)$ is the electrostatic potential at each 
spatial position $\vect{x}$ in the medium and time $t$.

\subsection{Constitutive relations}
\label{sec:constitutive_laws}

In this section, we provide the mathematical characterization of
the fluxes $\vect{j}_i$ and $\vect{j}_T$ and of the other 
model parameters in system~\eqref{eq:basic_system}.
To this purpose, we follow the classical 
references~\cite{Landau1984,degrootmazur1984} and, 
for more recent applications,~\cite{cheremisin2001,consiglieri2013},  
and assume that both current density $\vect{j}_i$ and 
energy flux density $\vect{j}_T$ can be expressed as the sum
of two contributions, namely, an electrochemical flux and
a thermal flux, so that:
\begin{subequations}\label{eq:fluxes}
\begin{align}
\vect{j}_i = \vect{j}_i^{ec} + \vect{j}_i^{th}
& \qquad i=1, \ldots, M \label{eq:flux_ect_i} \\[2mm]
\vect{j}_T =  \vect{j}_T^{ec} + \vect{j}_T^{th}. &
\label{eq:flux_ect_T}
\end{align}
\end{subequations}

\subsubsection{Electrical fluxes}\label{sec:electrical_model}
Let $\sigma_i$ denote the electrical conductivity
of species $N_i$ defined as
\begin{equation}\label{eq:conductivity}
\sigma_i = q |z_i| \mu^{el}_i N_i \qquad i=1, \ldots, M
\end{equation}
where $\mu^{el}_i$ is the electrical mobility of the $i$-th species.

The electrochemical flux associated with $N_i$ is~\cite{bardfalkner2001}:
\begin{subequations}\label{eq:ec_flux}
\begin{align}
\vect{j}_i^{ec} = -\sigma_i \nabla \varphi_i^{ec} 
& \qquad i=1, \ldots, M \label{eq:ec_flux_i} \\[2mm]
\varphi_i^{ec} = \varphi + \Frac{\mu^c_i}{z_i F}  
& \qquad i=1, \ldots, M \label{eq:pot_ec_i}
\end{align}
\end{subequations}
where $\varphi_i^{ec}$ is the \emph{electrochemical potential} 
of the $i$-th species given by the sum of the electrical potential 
$\varphi$ and of the \emph{chemical potential} 
\begin{equation}\label{eq:chem_pot}
\varphi_i^c:= \Frac{\mu^c_i}{z_i F},
\end{equation}
$\mu_i^c$ and $F$ being the chemical energy of the $i$-th species
and Faraday's constant, respectively.
In a homogeneous material ($N_i=const$), $\mu_i^c$ is constant 
so that the electrochemical potential is just a constant shift of the 
electric potential. In a non-homogeneous material ($N_i \neq const$),
the chemical energy is defined as
\begin{equation}\label{eq:chem_energy}
\mu^c_i = R T \ln \left( \Frac{N_i}{N_{ref}} \right)
\qquad i=1, \ldots, M
\end{equation}
where $R$ is the ideal gas constant 
and $N_{ref}$ is a reference concentration, so that 
the electrochemical potential is~\cite{bardfalkner2001}
\begin{equation}\label{eq:ec_pot_non_homog}
\varphi_i^{ec} = \varphi + \Frac{R T}{z_i F} 
\ln \left( \Frac{N_i}{N_{ref}} \right) = 
\varphi + \Frac{K_B T}{z_i q} \ln \left( \Frac{N_i}{N_{ref}} \right)
\qquad i=1, \ldots, M
\end{equation}
where $K_B$ is Boltzmann's constant.

Let us now consider the thermal current density $\vect{j}_i^{th}$. 
We have~\cite{cheremisin2001}:
\begin{subequations}\label{eq:th_flux}
\begin{align}
\vect{j}_i^{th} = -\sigma_i \nabla \varphi^{th} 
& \qquad i=1, \ldots, M \label{eq:th_flux_i} \\[2mm]
\varphi^{th} = \alpha T \label{eq:pot_th}
\end{align}
\end{subequations}
where $\varphi^{th}$ is the \emph{thermal potential}, 
$\alpha$ being the thermopower coefficient of the material.

Gathering together the above definitions of the
various flux and potential contributions, we can write
a \emph{generalized Ohm's law} for the current density
associated with the $i$-th chemical:
\begin{subequations}\label{eq:curr_density_ohm}
\begin{align}
\vect{j}_i = \sigma_i \vect{E}_{i}^{thec}
& \qquad i=1, \ldots, M \label{eq:curr_i} \\[2mm]
\vect{E}_{i}^{thec}:= - \nabla \psi_i & \qquad i=1, \ldots, M 
\label{eq:E_th_ec_i} \\[2mm]
\psi_i:= \varphi_i^{ec} + \varphi_i^{th} =  
\varphi + \Frac{\mu^c_i}{z_i F} + \alpha T
& \qquad i=1, \ldots, M \label{eq:pot_tot_i}
\end{align}
\end{subequations}
where $\psi_i$ is the \emph{thermo-electrochemical} potential
of the $i$-th chemical and $\vect{E}_{i}^{thec}$ is the 
\emph{thermo-electrochemical field} experienced by the $i$-th chemical. 
\begin{remark}[The generalized Drift-Diffusion model]
\label{rem:generalized_DD}

Replacing~\eqref{eq:conductivity} into~\eqref{eq:curr_density_ohm}
we obtain
the following equivalent form of the current density 
associated with the $i$-th chemical:
\begin{subequations}\label{eq:curr_density_DD}
\begin{align}
\vect{j}_i = q |z_i| \mu_i^{el} N_i \vect{E}_i^{el} 
- q z_i D_i \nabla N_i 
& \qquad i=1, \ldots, M \label{eq:curr_i_DD} \\[2mm]
\vect{E}_{i}^{el}:= \vect{E} - \alpha \nabla T -
\varphi_i^{c} \Frac{\nabla T}{T} & \qquad i=1, \ldots, M 
\label{eq:E_el_eff_i} \\[2mm]
D_i:= \Frac{K_B T}{q |z_i|} \mu_i^{el} 
& \qquad i=1, \ldots, M \label{eq:generalized_einstein}
\end{align}
\end{subequations}
where $D_i$ is the \emph{generalized diffusion coefficient}
of the $i$-th chemical, related to the electrical mobility $\mu_i^{el}$
through the \emph{generalized Einstein 
relation}~\eqref{eq:generalized_einstein} and 
$\vect{E}_i^{el}$ is the \emph{generalized electric field}
experienced by the chemical $N_i$.
Thus, Eq.~\eqref{eq:curr_i_DD} represents the 
\emph{generalized Drift-Diffusion (DD) model} for ionic charge 
transport in a non-homogeneous and non-isothermal material.
If the material is in isothermal conditions and 
electrons and holes are considered for transport, 
relation~\eqref{eq:curr_i_DD} 
degenerates into the classical DD model~\cite{Jerome1996}. 
In this case, two chemicals are flowing in the material
($M=2$), namely, negatively charged electrons ($z_1=-1$)
and positively charged holes ($z_2=+1$). 
\end{remark}

\subsubsection{Thermal flux}\label{sec:thermal_model}
Let $\kappa$ denote the thermal conductivity of 
the material. Then, classical Fourier law states that
heat thermal flow in the material is expressed by
the following relation
\begin{equation}\label{eq:heat_flux}
\vect{j}_T^{th} = -\kappa \nabla T.
\end{equation}
Heat is also transported in the direction of the total current 
flow in the material according to the following 
relation~\cite{degrootmazur1984}

\begin{equation}\label{eq:convective_heat_flow}
\vect{j}_T^{ec} = \psi \vect{j}
\end{equation}
where:
\begin{subequations}\label{eq:thermal_field_parameter}
\begin{align}
\psi:= \varphi + \alpha T + \Sum{i=1}{M} 
\Frac{K_B T}{z_i q} \ln \left( \Frac{N_i}{N_{ref}} \right)
& \label{eq:total_th_ec_pot} \\[2mm]
\vect{j}:= \Sum{i=1}{M} \, \vect{j}_{i}& 
\label{eq:total_curr_i} 
\end{align}
\end{subequations}
are the \emph{total} thermo-electrochemical potential and 
current density, respectively.
Gathering together the above definitions of the
various flux and potential contributions, we can write
the thermo-electrochemical heat flux in a concise
advection-diffusion form
\begin{equation}\label{eq:heat_thec}
\vect{j}_T = \psi \vect{j}  - \kappa \nabla T.
\end{equation}

\subsection{Model coefficients, sources and sinks}
\label{sec:coeff_sources_sinks}

To complete the description of the thermo-electrochemical
model we need specify the mathematical form of 
the physical parameters and coefficients. 
For sake of simplicity, we assume henceforth that 
the net production rates $R_i$ and $\mathcal{Q}_T$ are identically
equal to zero. These assumptions are equivalent to state that
sources and sinks in the material bulk accounting for mutual 
interactions among the chemicals are neglected in our description.
Concerning the other model parameters, we assume from now on 
that the electrical mobilities $\mu_i^{el}$, the 
thermopower $\alpha$, the dielectric permittivity $\varepsilon$, 
the thermal conductivity $\kappa$, the mass density and the
specific heat $\rho$ and $c$ are constant positive quantities
whose numerical values are specified in 
Sect.~\ref{sec:simulations}.

\begin{example}[The case of silicon devices]
A significant example of the application of the 
thermo-electrochemical model illustrated in this section is
provided by the study of silicon devices traditionally employed
in the semiconductor technology for microelectronics applications.
The corresponding version of system~\eqref{eq:basic_system}
including Joule heat dissipation but not thermo-electric power 
effects, is usually referred 
to as \emph{Energy-Transport} (ET) model (see~\cite{Jerome1996}).
The extension of the ET to cover also thermo-electric mechanisms
(Peltier and Thomson effects) can be found in~\cite{Baccarani1993}.
\end{example}

\section{Geometry and structure of the device}\label{sec:geometry}

\begin{figure}[h!]
\centering
\subfigure[3D view]
{\includegraphics[width=0.35\textwidth,height=0.55\textwidth]{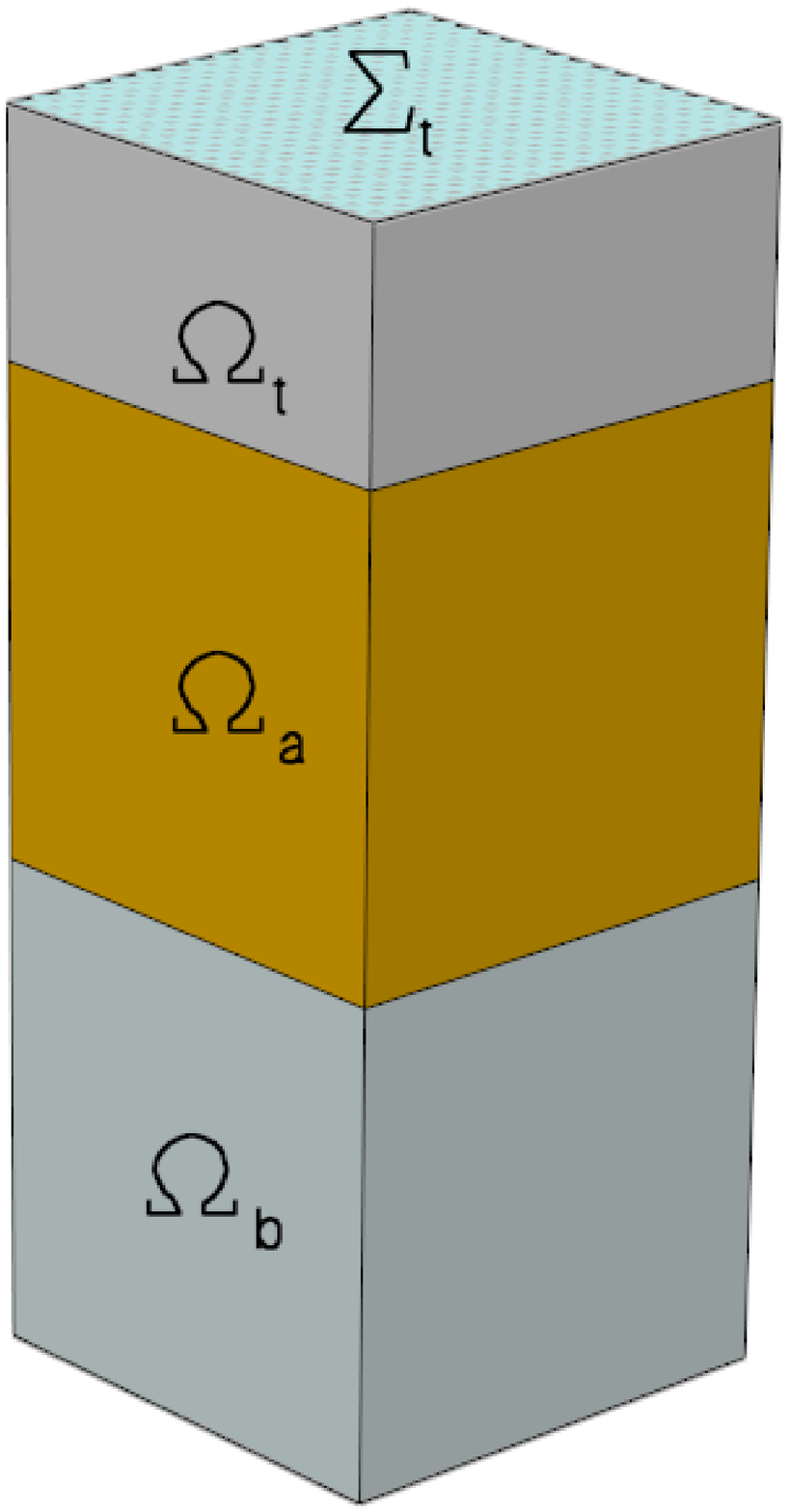}
\label{fig:cube_3d}}
\subfigure[2D view]
{\includegraphics[width=0.45\textwidth,height=0.55\textwidth]{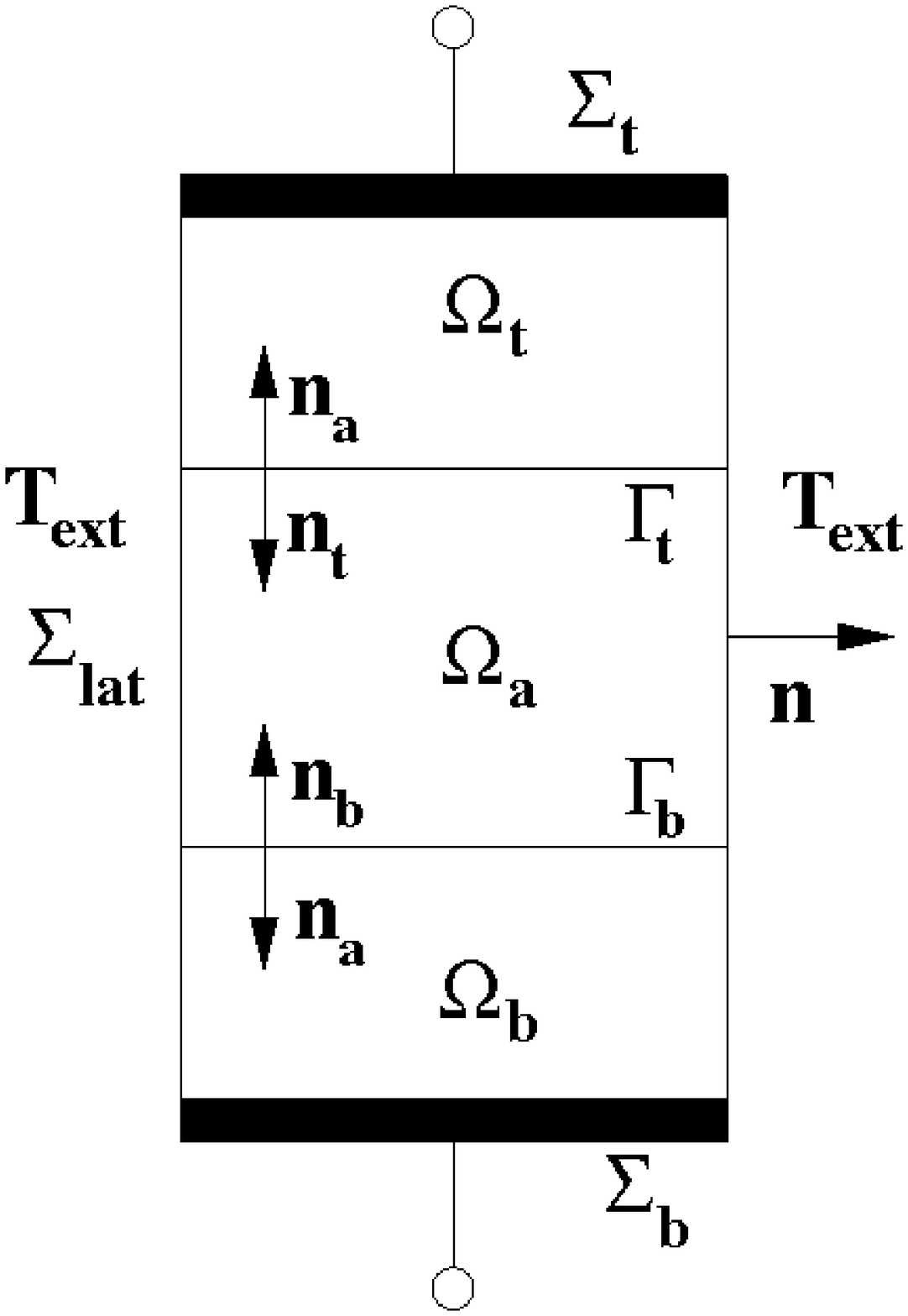}
\label{fig:geo_2d}}
\caption{Geometry of a typical semiconductor device for electronics 
applications. Left: three-dimensional scheme. Right: 
simplified two-dimensional scheme.}
\label{fig:geo_3d}
\end{figure}

In this section we address the geometrical description of
the semiconductor device object of the present work.
Fig.~\ref{fig:cube_3d} shows a perspective view of a typical 3D 
template devices for electronics applications. 
More complex device configurations will be investigated in
Sect.~\ref{sec:simulations}.
The device is characterized by an intrinsically material heterogeneous
structure composed of 
an active region (yellow layer) sandwiched between two 
inactive regions (red and grey blocks). These latter regions 
accomplish several important functions: 
1) they provide electrical and thermal
connection with the external environment, allowing to 
apply a voltage and thermal drop across the device; 
2) one of them supplies the intermediate active region 
with the appropriate thermo-electrochemical driving energy; 
and 3) the other one collects the thermo-electrochemical
current flux produced by the active region and transfers it to
the external circuit connected in series to the device
for further use.

In view of the mathematical modeling of the problem, we
consider in the present article the slightly simplified 
schematic geometrical representation illustrated, with a
two-dimensional (2D) cut view, in Fig.~\ref{fig:geo_2d}.
The device region is an open bounded domain $\Omega \subset 
\mathbb{R}^3$ consisting of the union of three subdomains:
the two inactive regions, $\Omega_t$ (top) and 
$\Omega_b$ (bottom), 
and the intermediate active layer, $\Omega_{a}$. 
The external boundary of the device, $\partial \Omega$, is made of
the union of three open disjoint surfaces, $\Sigma_t$, $\Sigma_b$ and
$\Sigma_{lat}$, on which an outward unit normal vector $\vect{n}$
is defined.

The top and bottom
surfaces, $\Sigma_t$ and $\Sigma_b$, are the electrical 
and thermal contacts
where external voltage and thermal sources are applied.
The lateral surface of the device, $\Sigma_{lat}$, is the
material interface between the device and the external environment
whose temperature is denoted by $T_{ext}$. The portion
of $\Sigma_{lat}$ belonging to the boundary of the active region
is denoted henceforth as $\Sigma_{lat}^a$.

The intermediate subdomain $\Omega_{a}$ 
is separated from the top inactive region by the interface surface 
$\Gamma_t$ on which we define the normal unit vectors $\vect{n}_t$ 
and $\vect{n}_{a}$, outwardly directed 
on each surface with respect to $\Omega_t$ and 
$\Omega_{a}$. 
In a similar manner, on the bottom interface surface $\Gamma_b$
separating $\Omega_{a}$ from the inactive region $\Omega_b$
we define the normal unit vectors $\vect{n}_b$ 
and $\vect{n}_{a}$, outwardly directed 
on each surface with respect to $\Omega_b$ and $\Omega_{a}$.
Clearly, $\vect{n}_t = - \vect{n}_{a}$ on $\Gamma_t$ and 
$\vect{n}_b = - \vect{n}_{a}$ on $\Gamma_b$.

\section{Mathematical Model of a 3D Semiconductor Device}
\label{sec:model}

In this section we use the basic theory developed in
Sect.~\ref{sec:th_el_ch} to construct the thermo-electrochemical 
mathematical model of the 3D heterogeneous semiconductor device 
introduced in Sect.~\ref{sec:geometry}. 
From now on, we denote by $(0,t_{fin})$ the time interval in which we 
study the dynamical behaviour of the device and
we assume that:
\begin{description}
\item[(A1)] the two inactive regions $\Omega_b$ 
and $\Omega_t$ are metals or degenerate semiconductors;
\item[(A2)] charge transport in the device active region 
$\Omega_{a}$ is only due to electrons that are injected by tunneling 
from one of the inactive regions into the active layer.
\end{description}
Based on (A2), we limit ourselves to considering the case 
$M=1$ and set $N_1:= n$, $n=n(\vect{x},t)$ being the 
electron number density in the device. Based on (A1), we also set:
\begin{equation}\label{eq:conc_metals}
n(\vect{x},t)=
\left\{
\begin{array}{ll}
\overline{n}^b & \qquad \forall \vect{x} \in \Omega_b 
\quad \forall t >0 \\[2mm]
\overline{n}^t & \qquad \forall \vect{x} \in \Omega_t
\quad \forall t >0 
\end{array}
\right.
\end{equation}
where $\overline{n}^b$ and $\overline{n}^t$ are the constant values 
of the electron concentration in the two metallic domains.
Relations~\eqref{eq:conc_metals} express the fact that 
the two metal regions $\Omega_b$ and $\Omega_t$ 
are homogeneous conductors, so that electrical conductivity 
is uniform (equal to $\sigma^b$ and $\sigma^t$, respectively)
and electric current transport is governed by the ideal Ohm's law.
Finally, we assume that the chemical energies 
of the bottom and top metal regions are constant values 
equal to $\overline{\mu}_c^b$ and 
$\overline{\mu}_c^t$, respectively, the thermopower 
coefficient $\alpha$ is a piecewise constant function equal to
$\alpha^b$ and $\alpha^t$ in the bottom and 
metal regions, respectively,
and to $\alpha^a$ in the active region, the electron 
electrical mobility $\mu_n^{el}$ is a positive constant in the
active region and the mass density $\rho$ and specific heat
$c$ are piecewise constant positive functions equal to 
$\rho^b$, $\rho^a$, $\rho^t$ and $c^b$, $c^a$, $c^t$,  
respectively.

Replacing~\eqref{eq:curr_density_DD} and~\eqref{eq:heat_thec} 
into~\eqref{eq:basic_system} and using~\eqref{eq:conc_metals}, 
we end up with the following system of PDEs in conservation
form to describe the thermo-electrochemical processes occurring
in a 3D heterogeneous semiconductor device:
\begin{subequations}\label{eq:model_thec}
\begin{align}
\Div \vect{j}_{\varphi} = f_{\varphi}
& \qquad \mbox{in } \Omega \times (0, t_{fin}) 
\label{eq:Poisson_generalized} \\[2mm]
- q \Frac{\partial n}{\partial t} + \Div \vect{j}_n = 0 
& \qquad \mbox{in } \Omega_{a} \times (0, t_{fin}) 
\label{eq:electron_continuity} \\[2mm]
\Frac{\partial}{\partial t}  (\rho c T) + 
\Div \vect{j}_T = 0
& \qquad \mbox{in } \Omega \times (0, t_{fin}).
\label{eq:heat_generalized}
\end{align}
\end{subequations}

The chosen ordering of the equations in system~\eqref{eq:model_thec} 
reflects the structure of the solution map that is used to
iteratively solve the problem as illustrated in
Sect.~\ref{sec:comp_tech}.

Eq.~\eqref{eq:Poisson_generalized} represents a 
\emph{generalized} Poisson equation in the whole 
device, because it coincides with the 
differential Gauss' law~\eqref{eq:Poisson} in the active region 
while in the two metal regions it takes the form and meaning of 
an electro-thermal Ohm's law, having defined in a piecewise 
manner over the device domain $\Omega$ 
the electro-thermal flux $\vect{j}_\varphi$ as:
\begin{subequations}\label{eq:fluxes_ohmic_diel_thec}
\begin{align}
\vect{j}_\varphi = 
\left\{
\begin{array}{ll}
-\sigma^b \nabla \varphi - \alpha^b \sigma^b \nabla T & \qquad
\mbox{in } \Omega_b \times (0, t_{fin}) \\[2mm]
-\varepsilon^{a} \nabla \varphi & \qquad
\mbox{in } \Omega_{a} \times (0, t_{fin}) \\[2mm]
-\sigma^t \nabla \varphi - \alpha^t \sigma^t \nabla T& \qquad
\mbox{in } \Omega_t \times (0, t_{fin}) \\[2mm]
\end{array}
\right. 
& \label{eq:electrothermal_flux_thec}
\end{align}
and the space charge density $f_\varphi$ as:
\begin{align}
f_\varphi = 
\left\{
\begin{array}{ll}
0 & \qquad
\mbox{in } \Omega_b \times (0, t_{fin}) \\[2mm]
-q n + q \mathcal{D} & \qquad
\mbox{in } \Omega_{a} \times (0, t_{fin}) \\[2mm]
0 & \qquad
\mbox{in } \Omega_t \times (0, t_{fin}).
\end{array}
\right. 
& \label{eq:spacechargedensity_thec}
\end{align}
\end{subequations}
Eq.~\eqref{eq:electron_continuity} is the electron continuity
equation in the active region, the electron current
density $\vect{j}_n$ being defined as
\begin{equation}\label{eq:electron_flux_thec}
\vect{j}_n= -q \mu_n^{el} n \left[\nabla \varphi 
+ \alpha^{a} \nabla T - \Frac{K_B T}{q} 
\ln \left( \Frac{n}{N_{ref}} \right) \Frac{\nabla T}{T} \right]
+ q D_n \nabla n.
\end{equation}
Eq.~\eqref{eq:heat_generalized} is the heat flow equation in the whole 
device structure, the total heat flux $\vect{j}_T$ being defined as
\begin{subequations}\label{eq:thermal_fluxes_thec}
\begin{align}
\vect{j}_T = \psi_n \vect{j} - \kappa \nabla T 
& \qquad \mbox{in } \Omega \times (0, t_{fin}) 
\label{eq:heat_flux_thec}
\end{align}
where the thermo-electrochemical potential $\psi_n$
and density flux $\vect{j}$
are defined in a piecewise manner over the device domain $\Omega$ as:
\begin{equation}\label{eq:psi_thec}
\psi_{n} = 
\left\{
\begin{array}{ll}
\varphi - \Frac{\overline{\mu}^b_c}{F} + \alpha^b T & \qquad
\mbox{in } \Omega_b \times (0, t_{fin}) \\[4mm]
\varphi - \Frac{K_B T}{q} \log \left( \Frac{n}{N_{ref}} \right)
+ \alpha^a T & \qquad
\mbox{in } \Omega_{a} \times (0, t_{fin}) \\[2mm]
\varphi - \Frac{\overline{\mu}^t_c}{F} + \alpha^t T & \qquad
\mbox{in } \Omega_t \times (0, t_{fin})
\end{array}
\right.
\end{equation}
and:
\begin{equation}\label{eq:thermoelectrical_flux_thec}
\vect{j} = 
\left\{
\begin{array}{ll}
-\sigma^b \nabla (\varphi + \alpha^b T) & \qquad
\mbox{in } \Omega_b \times (0, t_{fin}) \\[2mm]
\vect{j}_n & \qquad \mbox{in } \Omega_{a} \times (0, t_{fin}) \\[2mm]
-\sigma^t \nabla (\varphi + \alpha^t T) & \qquad
\mbox{in } \Omega_t \times (0, t_{fin}).
\end{array}
\right. 
\end{equation}
\end{subequations}

To complete the mathematical model of thermo-electrochemical 
transport in a semiconductor device, we need specify suitable initial
and boundary conditions. 

Concerning the initial conditions, we set:
\begin{subequations}\label{eq:ICs}
\begin{align}
n(\vect{x}, 0) = n_0(\vect{x}) & \qquad \forall \vect{x} \in 
\Omega_{a} \label{eq:IC_n} \\[2mm]
T(\vect{x}, 0) = T_0(\vect{x}) & \qquad \forall \vect{x} \in 
\Omega \label{eq:IC_T}
\end{align}
\end{subequations}
where $n_0:\Omega_{a} \rightarrow \mathbb{R}$ and 
$T_0: \Omega \rightarrow \mathbb{R}$ are given positive functions.

Concerning the boundary conditions, for all $t \in (0, t_{fin})$ we set:
\begin{subequations}\label{eq:BCs}
\begin{align}
\varphi = \overline{\varphi}^b & \qquad \mbox{on } \Sigma_b
\label{eq:BC_D_phi_b} \\[2mm]
\varphi = \overline{\varphi}^t & \qquad \mbox{on } \Sigma_t
\label{eq:BC_D_phi_t} \\[2mm]
\vect{j}_\varphi \cdot \vect{n} = 0 
& \qquad \mbox{on } \Sigma_{lat} \label{eq:BC_N_J_phi}
\end{align}
for the generalized Poisson equation, and:
\begin{align}
-\vect{j}_n \cdot \vect{n}_b = j_{tunnel}  
& \qquad \mbox{on } \Gamma_{b}  
\label{eq:BC_N_n_b} \\[2mm]
-\vect{j}_n \cdot \vect{n} = q v_{eq} (n-n_{eq})
& \qquad \mbox{on } \Gamma_{t} 
\label{eq:BC_N_n_t} \\[2mm]
\vect{j}_n \cdot \vect{n} = 0 
& \qquad \mbox{on } \Sigma_{lat}^a
\label{eq:BC_N_n_lat}
\end{align}
for the electron continuity equation, and:
\begin{align}
T = \overline{T}^b & \qquad \mbox{on } \Sigma_b
\label{eq:BC_D_T_b} \\[2mm]
T = \overline{T}^t & \qquad \mbox{on } \Sigma_t
\label{eq:BC_D_T_t} \\[2mm]
\vect{j}_T \cdot \vect{n} = \gamma_T (T - T_{ext})
& \qquad \mbox{on } \Sigma_{lat} 
\label{eq:BC_N_T_lat}
\end{align}
\end{subequations}
for the generalized heat equation.

Let us address the mathematical and physical interpretation
of the above boundary conditions.

Relations~\eqref{eq:BC_D_phi_b}-~\eqref{eq:BC_D_phi_t} are
non-homogeneous Dirichlet boundary conditions for the electric
potential expressing the physical fact that
the electric contacts are equipotential surfaces equal 
to the externally applied voltage sources 
$\overline{\varphi}^b$ and $\overline{\varphi}^t$.

Relations~\eqref{eq:BC_N_J_phi}-~\eqref{eq:BC_N_n_lat}
are homogeneous Neumann conditions expressing the physical
fact that charge transport in the device is self-contained,
i.e., current lines start and close between the
two bottom and top surfaces.

Relation~\eqref{eq:BC_N_n_b} is a non-homogeneous 
Neumann condition for the electron flux, $j_{tunnel}$ being
the (positive) electron current density injected by tunneling from
the bottom metal region $\Omega_b$ into the active layer
$\Omega_{a}$ across the separating surface $\Gamma_b$. 

Relation~\eqref{eq:BC_N_n_t} is a Robin boundary condition
expressing the net electron current flux flowing  
between the active region and the top metal region
across the separating interface surface $\Gamma_t$.
The mathematical form of this boundary condition is analogous
to that used to describe current flux balance 
at a Schottky interface between a metal and a 
semiconductor accounting for thermoionic emission (from the metal) 
and drift-diffusion injection (from the semiconductor)~\cite{sze2006}. 
According to this interpretation, $v_{eq}$ and $n_{eq}$ are
the values of drift velocity and electron concentration at thermodynamical
equilibrium conditions while $n$ is the unknown value of electron number 
density on the interface side of the active layer.

Relations~\eqref{eq:BC_D_T_b}-~\eqref{eq:BC_D_T_t} are
non-homogeneous Dirichlet boundary conditions for device 
temperature expressing the physical fact that
the electric contacts are also equi-thermal surfaces equal 
to the externally applied positive thermal sources $\overline{T}^b$ and
$\overline{T}^t$.

Relation~\eqref{eq:BC_N_T_lat}
is a Robin boundary condition expressing the net heat flux 
exchange between the device and the surrounding environment,
$\gamma_T$ being a non-negative heat transfer coefficient.
 
\begin{remark}[The nature of the model]
It is interesting to notice that the coupled set of 
equations~\eqref{eq:model_thec}-~\eqref{eq:fluxes_ohmic_diel_thec} 
supplemented by the initial conditions~\eqref{eq:ICs} 
and boundary conditions~\eqref{eq:BCs} constitute an incompletely
parabolic system of PDEs because of the need of satisfying
the elliptic constraint~\eqref{eq:Poisson_generalized} at each time level.
This issue makes the treatment of the problem quite difficult,
both in analytical and numerical terms.
\end{remark}

\section{Computational Techniques}\label{sec:comp_tech}

In this section we describe the various steps that
transform the PDE model~\eqref{eq:model_thec}-~\eqref{eq:BCs}
into the successive solution of linear systems of algebraic 
equations of large size that represent the discrete counterpart
of the problem.

\subsection{Time semi-discretization}\label{sec:time_discretization}

We divide the time interval $(0, t_{fin})$ into 
a finite number $N_t \geq 1$ of time slabs of nonuniform
width $\Delta t_k := t_k - t_{k-1}$, $k=1, \ldots, N_t$ with $t_0=0$, 
in such a way that discrete time levels are denoted as 
$t_k$, $k=0, \ldots, N_t$.
The choice of a nonuniform discretization of the time variable 
is made in order to properly track the wide dynamical range 
of the temporal scales of the thermo-electrochemical phenomena
occuring in the device under investigation which may vary between
nanoseconds to milliseconds up to even seconds. 
In the present computer implementation, the sequence of values of 
$\Delta t_k$ is user-defined and for time advancing the Backward Euler 
(BEM) method is adopted because of its unconditional 
stability. An alternative approach 
based on the use of higher-order methods coupled with 
adaptive strategies for automatic time-step selection 
(see, e.g.,~\cite{bank1985,Lambert1991,ascher1998cmo}) will be considered 
in a future extension of the computational scheme proposed in the 
present article. 

\subsection{Solution map}\label{sec:gummel_map}
Throughout the remainder of the article, given a function 
$f=f(\vect{x},t)$ we set $f_k(\vect{x}):=f(\vect{x}, t_k)$ 
for every $k=0, \ldots, N_t$. We also denote by $\chi^b$, $\chi^a$
and $\chi^t$ the characteristic functions of the sets $\Omega^b$,
$\Omega^a$ and $\Omega^t$, respectively, 
such that $\chi^\nu(\vect{x}) =1$ if 
$\vect{x} \in \Omega_\nu$ and $\chi^\nu(\vect{x}) =0$ if 
$\vect{x} \notin \Omega_\nu$, $\nu=\left\{b, a, t\right\}$. 
The functional iteration illustrated
below is used to linearize the thermo-electrochemical model 
upon previous time 
semidiscretization with the BEM: \\[2mm]
given $[n_k, T_k]^t$, $k=0, \ldots, N_t-1$, execute the following
solution steps: \\

\noindent
A) set $n^{(0)}:=n_k, \, T^{(0)}:=T_k$;

\noindent
B) for $m=0, 1, \ldots,$ until convergence, solve: \\
\begin{subequations}\label{eq:gummel_thec}
\begin{align}
\Div \left(-\mathcal{S}^\varphi \nabla U \right) 
= \left(-q n^{(m)} + q \mathcal{D}\right) \chi^a
-\Div \left(-\mathcal{S}^T \nabla T^{(m)} \right)
& \qquad \mbox{in } \Omega 
\label{eq:Poisson_generalized_m}
\end{align}
with:
\begin{align}
\mathcal{S}^\varphi:= \sigma^b \chi^b + \varepsilon^{a} \chi^{a} 
+ \sigma^t \chi^t & \qquad \mbox{in } \Omega \nonumber \\[2mm]
\mathcal{S}^T:= \sigma^b \alpha^b \chi^b 
+ 0 \cdot \chi^{a} + \sigma^t \alpha^t \chi^t & 
\qquad \mbox{in } \Omega \nonumber 
\end{align}
and set $\varphi^{(m+1)}:=U$; 
\begin{align}
\Frac{q U}{\Delta t_k} 
+ \Div \left( \vect{V}_n^{(m+\frac{1}{2})} U - q D_n \nabla U \right)
= \Frac{q n^{(m)}}{\Delta t_k} & \qquad \mbox{in } \Omega_{act}
\label{eq:electron_continuity_m}
\end{align}
with:
\begin{align}
\vect{V}_n^{(m+\frac{1}{2})}:= q \mu_n^{el} \left[\nabla \varphi^{(m+1)} 
+ \alpha^{a} \nabla T^{(m)} - \Frac{K_B T^{(m)}}{q} 
\ln \left( \Frac{n^{(m)}}{N_{ref}} \right) 
\Frac{\nabla T^{(m)}}{T^{(m)}} \right] & \nonumber
\end{align}
and set $n^{(m+1)}:=U$; 
\begin{align}
\Frac{\rho c U}{\Delta t_k} +
\Div \left(\vect{V}_T^{(m+\frac{1}{2})} U  - \kappa \nabla U \right)
= \Frac{\rho c T^{(m)}}{\Delta t_k} 
-\Div \left(\varphi_{ec}^{(m+\frac{1}{2})} \vect{V}_T^{(m+\frac{1}{2})} 
\right) & \qquad \mbox{in } \Omega.
\label{eq:heat_generalized_m}
\end{align}
with:
\begin{align}
\vect{V}_T^{(m+\frac{1}{2})}:=
\alpha^b \left[ -\sigma^b \nabla (\varphi^{(m+1)} 
+ \alpha^b T^{(m)}) \right] \chi^b 
+ \alpha^a \vect{j}_n(\varphi^{(m+1)}, n^{(m+1)}, T^{(m)}) 
\chi^a & \nonumber \\[2mm]
+ \alpha^t \left[ -\sigma^t \nabla (\varphi^{(m+1)} 
+ \alpha^t T^{(m)}) \right] \chi^t
& 
\qquad \mbox{in } \Omega  \nonumber \\[2mm]
\varphi_{ec}^{(m+\frac{1}{2})}:= \left( \varphi^{(m+1)} 
- \Frac{\overline{\mu}_c^b}{F} \right) \chi^b 
+  \left[ \varphi^{(m+1)} -\Frac{K_B T^{(m)}}{q}
\ln \left( \Frac{n^{(m+1)}}{N_{ref}}\right) \right] 
\chi^a & \nonumber \\[2mm] 
+ \left( \varphi^{(m+1)} 
- \Frac{\overline{\mu}_c^t}{F} \right) \chi^t
& \qquad \mbox{in } \Omega  \nonumber
\end{align}
and set $T^{(m+1)}:=U$. 

\null
\noindent
C) Let $\vect{U}:= \left[ \varphi, n, T \right]^t$ denote the 
solution triple. Should the sequence $\left\{ \vect{U}^{(m)} \right\}$
be converging to a fixed point $\vect{U}^{\ast}$, then set:
\begin{equation}\label{eq:update_time_solution}
\varphi_{k+1}:= \varphi^{\ast}, \quad n_{k+1}:= n^{\ast}, \quad
T_{k+1}:= T^{\ast}
\end{equation}
and proceed to the next time level.
\end{subequations}

\null

\begin{figure}[h!]
\centering
\includegraphics[width=0.45\textwidth,height=0.85\textwidth]{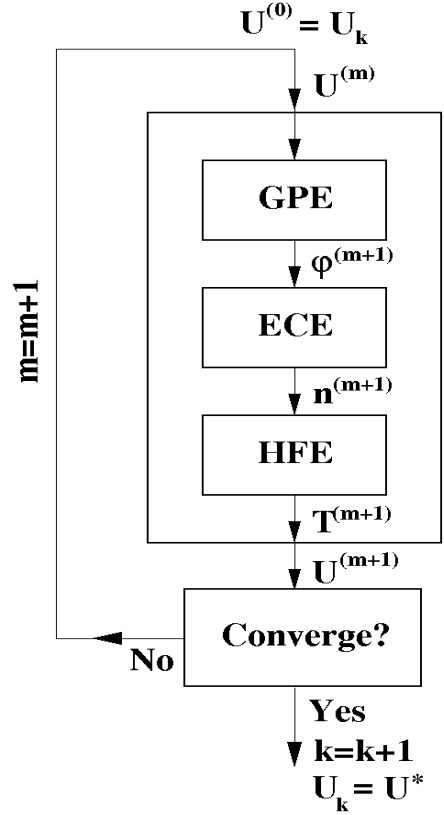}
\caption{Flow-chart of a single step of the solution map. If 
the iteration has reached convergence to a fixed point $\vect{U}^\ast$, 
then the solution vector at next time level is set equal to $\vect{U}^\ast$
and the algorithm advances to the next time level $t_{k+1}$.}
\label{fig:gummel_map}
\end{figure}

The solution map~\eqref{eq:gummel_thec} can be regarded as 
the extension of the Gummel decoupled iteration 
widely employed in contemporary semiconductor device 
simulation (see~\cite{Selberherr84}) and thoroughly 
analyzed in~\cite{Jerome1996}. The study of existence (and possible
uniqueness) of a fixed point $\vect{U}^{\ast}$ 
of~\eqref{eq:gummel_thec} and of the convergence of the 
solution map as a function of model physical parameters
goes beyond the scope of the present article and will be
the object of a future publication.

The three steps A), B) and C) of the functional 
iteration~\eqref{eq:gummel_thec} are schematically
represented in the flow-chart of Fig.~\ref{fig:gummel_map} where
the symbol GPE in the first block indicates the
linear Generalized Poisson Equation~\eqref{eq:Poisson_generalized_m}, 
while ECE and HFE denote the linear Electron Continuity 
Equation~\eqref{eq:electron_continuity_m} and the
linear Heat Flow Equation~\eqref{eq:heat_generalized_m}, respectively.
The criterion adopted to monitor the convergence of the 
solution map~\eqref{eq:gummel_thec} is to stop the algorithm 
at the first value $m^{\ast} \geq 0$ of the iteration counter $m$ 
such that
$$
\| \vect{n}^{(m^{\ast}+1)} - \vect{n}^{(m^\ast)} \|_{2} < {\tt toll}
$$
where ${\tt toll}$ is a prescribed tolerance, 
$\vect{n}$ denotes the vector of nodal degrees of freedom 
of the finite element approximation $n_h$ and 
$$
\| \vect{w} \|_2:= \left(\Sum{i=1}{p} w_i^2\right)^{1/2}
$$
is the $2$-norm of a vector $\vect{w} \in \mathbb{R}^p$.
In the numerical experiments we have set ${\tt toll}= 10^{-3}$.

\subsection{Numerical approximation}\label{sec:FE_discretization}

In this section we carry out the numerical approximation 
of each linear boundary value problem in the Gummel iterative 
process~\eqref{eq:gummel_thec}
using the Galerkin Finite Element (FE) method. 
To this purpose, we introduce a partition $\mathcal{T}_h$ 
of the domain $\Omega$
into regularly shaped~\cite{Ciarlet1978} tetrahedral elements $K$
of average size $h$, $h >0$ denoting the discretization parameter.
On the triangulation $\mathcal{T}_h$, we define the 
finite dimensional space
\begin{equation}\label{eq:V_h}
V_h := \left\{v \in C^0(\overline{\Omega}) \, | \, 
v_h \in \mathbb{P}_1(K) \, \mbox{for all } K \in \mathcal{T}_h \right\}
\end{equation}
of piecewise affine functions that are continuous 
over the computational domain. The dimension of $V_h$ is denoted
henceforth by $N_h$ and coincides with the number of vertices of
$\mathcal{T}_h$. 

The standard Galerkin FE method consists of finding the 
approximation $U_h \in V_h$ of the weak solution $U$ 
of each problem in~\eqref{eq:gummel_thec}
(see~\cite{QV1994}) and gives rise to the solution of a 
linear system of algebraic equations 
\begin{equation}\label{eq:linear_system_FE}
\mathbf{K} \mathbf{U} = \mathbf{F}
\end{equation}
where $\mathbf{K} \in \mathbb{R}^{N_h \times N_h}$ 
is the stiffness matrix, $\mathbf{U} \in \mathbb{R}^{N_h}$ is the vector 
of nodal values $U_i$, $i=1, \ldots, N_h$, while 
$\mathbf{F} \in \mathbb{R}^{N_h}$ is the load vector.
The formulation may suffer of unwanted 
instabilities in the case where reaction and/or convection terms 
dominate over the diffusion term. Such instabilities typically 
show up under the form of spurious oscillations in the computed 
numerical solution which, in extreme cases, 
may even give rise to negative values of $U_h$. This latter event
is physically not acceptable, should $U_h$ represent a number
density or a temperature.

The simple-minded remedy to overcome these problems is to reduce 
the mesh size $h$, at the price, however,
of a considerable increase of the computational effort which
may become overwhelming in 3D simulations.
An alternative approach consists of introducing into the 
FE formulation suitable stabilization terms as discussed 
in~\cite{Hughes82,Hughes92,Hughes97,QV1994}.
These stabilized FE methods prevent (or strongly limit)
the onset of spurious oscillations preserving at the same time 
the optimal convergence properties of the plain FE approximation
but are not able, in general, to ensure the computed 
solution to be non-negative.

Since in our application (and more in general, in 
ther\-mo-e\-lec\-tro\-che\-mi\-cal
models) the property of $U_h$ of being 
non-negative is critical because of the physical meaning of the 
unknown (temperature, number density), in this article
we adopt the {\em exponentially fitted} or
{\em edge-averaged} finite element scheme thoroughly 
discussed and analyzed
in~\cite{BankBurglerFichtnerSmith1990,BankCoughranCowsar1998,Gatti1998} 
(in two spatial dimensions) 
and in~\cite{XuZikatanov1999,LazarovZikatanov2005,CdFThesis2006} 
(also in three spatial dimensions). The method 
is a multi-dimensional extension of the
classical Schar\-fet\-ter-Gum\-mel difference scheme~\cite{sg1969}
and gives rise to the linear algebraic system 
\begin{equation}\label{eq:linear_system_FE_stab}
\mathbf{K}^{SG} \mathbf{U}^{SG} = \mathbf{F}
\end{equation}
where $\mathbf{K}^{SG} \in \mathbb{R}^{N_h \times N_h}$ 
is the stiffness matrix associated with the exponentially fitted
discretization, $\mathbf{U}^{SG} \in \mathbb{R}^{N_h}$ is the 
corresponding solution vector while the right-hand side is the 
same as in~\eqref{eq:linear_system_FE}. 
According to Lemma 5.1 of~\cite{XuZikatanov1999}, it can be shown that
$\mathbf{K}^{SG}$ is a M-matrix~\cite{Varga1962} 
under suitable conditions on the 
shape regularity of the triangulation $\mathcal{T}_h$. 
This property implies the following important result which 
expresses the well-posedness and monotonicity of the discrete problem.
\begin{proposition}
\label{prop:well-posedness_h}
The linear algebraic system~\eqref{eq:linear_system_FE_stab} is
uniquely solvable. Moreover, if $F_i \geq 0$ for all $i=1, \ldots, N_h$
then the solution $\mathbf{U}^{SG}$ of~\eqref{eq:linear_system_FE_stab}
satisfies the discrete monotonicity property
\begin{equation}\label{eq:maximum_principle_h}
U_i^{SG} \geq 0 \qquad \forall i=1, \ldots, N_h.
\end{equation}
\end{proposition}

\section{Simulations and Results}\label{sec:simulations}
\begin{figure}[h!]
\centering
\subfigure[Template of a cubic structure]
{\includegraphics[width=0.35\textwidth,height=0.45\textwidth]{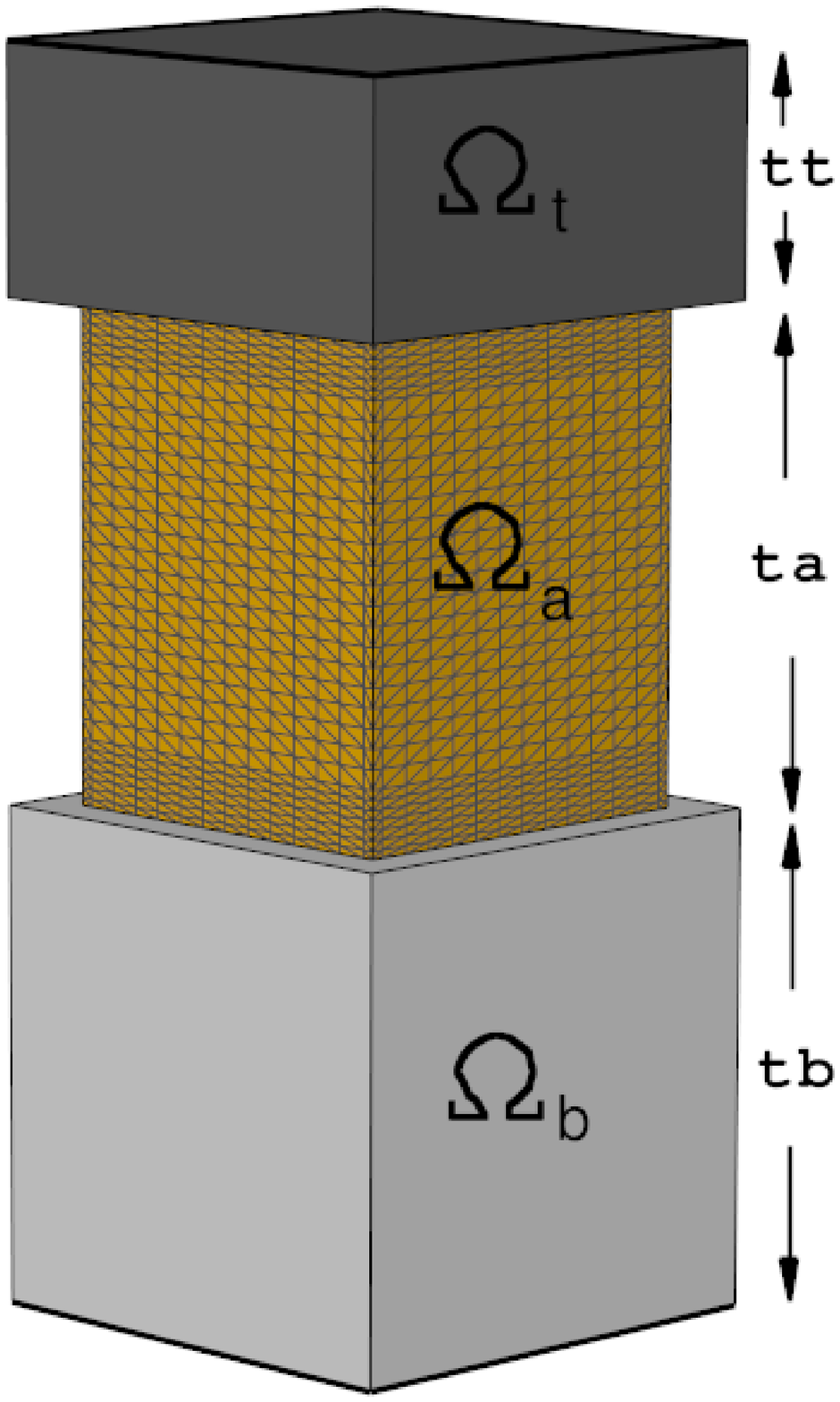}
\label{fig:cube_mesh}}
\subfigure[Template of a cylindrical structure]
{\includegraphics[width=0.45\textwidth,height=0.4\textwidth]{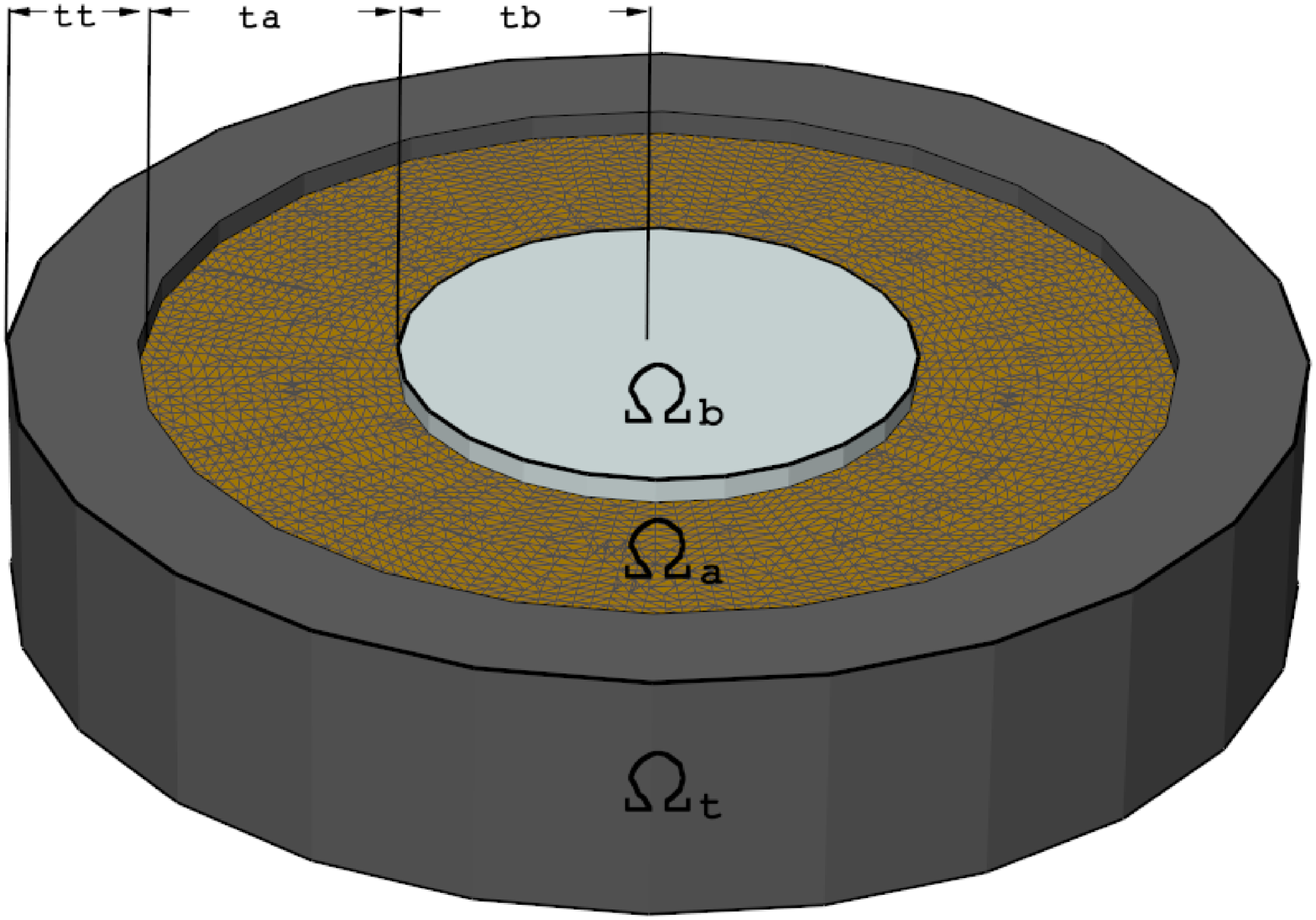}
\label{fig:cyl_mesh}}
\caption{Meshes of $\Omega_{a}$ used in the numerical experiments. 
Left: Cube with axis conformal mesh. 
Right: Cylinder with boundary conformal mesh.}
\label{fig:mesh}
\end{figure}

The model and the computational algorithm 
described in Sect.~\ref{sec:model} and Sect.~\ref{sec:comp_tech} 
have been implemented in a numerical code written in C++ 98
and compiled with gcc 4.5.2 in shared libraries on 64-bit architectures 
and run on multiple cores blades.
The code has then been applied to the simulation of 
several 3D structures with a cubic or cylindrical shape 
of which some examples are shown in Fig.~\ref{fig:mesh}.

Fig.~\ref{fig:mesh} visualizes typical meshes built with tetrahedra used 
in the discretization procedure.
Cubic structures are geometrically discretized 
using a general Delaunay mesh generator
because the most important surfaces are axis aligned.
In the case of cylindrical structures a Delaunay mesher 
able to build a boundary conformal mesh employing surface-adapted, 
anisotropic, mesh layers has been used in order to properly 
account for the influence of the 
interfaces on the boundary conditions~\cite{Synopsys}.
The resulting meshes are constituted by a number of elements 
varying from 90000 to 450000 depending on the simulated structure.
All the reported simulation results show the computed solutions
at steady-state conditions. In the case of $9 \, \unit{\mu m}$-size
structures the final time needed for steady-state to be reached is 
of the order of $100 \, \unit{\mu s}$, while in the case of cylindrical
$10 \, \unit{nm}$-size structures the final time needed for 
steady-state to be reached is of the order of $10\, \unit{\mu s}$.
The typical computational time for the presented cases  
varies between a minimum of 30 minutes to a couple of hours maximum.
In all reported simulation data and results, physical model parameters 
and variables are expressed in the units of the International System
according to the list of Sect.~\ref{sec:list}, except for the length 
scale which is expressed in $\unit{\mu m}$ for graphical convenience.

\subsection{Comparison with analytical solutions}\label{sec:analytical}
\begin{figure}[h!]
\centering
\includegraphics[width=0.65\textwidth,height=0.5\textwidth]
{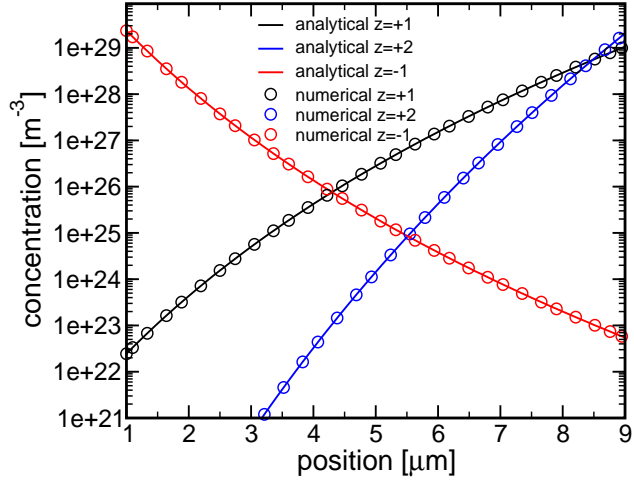}
\caption{Comparison between analytical and numerical 
solution (1D cut in the center of the $x-y$ plane along 
the $z$ axis) for Eq.~\eqref{eq:mass_balance_N_i}.
Three chemical species with charge $z=\pm1,+2$ are considered.}
\label{fig:chem_ana}
\end{figure}

This section is devoted to the comparison of the computed 
3D numerical solution with the 1D analytical solution obtained 
for simple cases. For this purpose, as reference structure we have 
used a cubic device characterized by different values of the 
thickness and imposed homogeneous Neumann boundary conditions 
on $\Sigma_{lat}$, $\Gamma_{t}$ and $\Gamma_{b}$.
Firstly, we have tested Eq.~\eqref{eq:mass_balance_N_i} with $R_{i}=0$ and 
a constant given electric field of strength equal to 
$|\textbf{E}|=10^6 \, \unit{V} \unit{m}^{-1}$ 
considering three different chemical species with charge 
$z=\pm q$ and $+2 q$ ($q$ being the electron charge).
For each species, the initial condition is set constant 
in all $\Omega_{a}$ and equal to $N_{i,0}=10^{28} \, \unit{m}^{-3}$.
A linear variation of temperature from $\Gamma_{t}=370 \, \unit{K}$ to 
$\Gamma_{b}=970 \, \unit{K}$ is imposed to the structure. 
Geometrical thicknesses are $t_{b}=t_{t}=0$ and $t_{a}=9 \unit{\mu m}$.
Fig.~\ref{fig:chem_ana} reports the results of the numerical simulation
(1D cuts along the $z$ axis in the center of the $x-y$ plane) compared 
with the analytical solution: symbols are for the numerical and lines 
for the analytical results.
No difference has been found between our implementation and the
exact stationary solution.

Secondly, we have considered Eq.~\eqref{eq:heat_generalized} with
an electric field directed along the $z$ axis, $\vect{E}=E_{0}\vect{z}$,
where $\vect{z}$ is the unit vector of the $z$ axis.
If an uniform concentration of electrons ($N_{e}$) is imposed in 
all $\Omega_{a}$ we can neglect the contribution coming from the 
diffusion term so that~\eqref{eq:thermal_fluxes_thec} reduces to
\begin{equation}\label{eq:heat_thec_simu}
\vect{j}_T = q \alpha N_{e} \mu_{e}^{el} E_{0}\vect{z} T  
- \kappa \nabla T
\end{equation}
where $\mu_{e}^{el}$ is the electron mobility. 
In this condition Eq.~\eqref{eq:heat_generalized}
can be solved analytically. The considered thicknesses of the cubic 
structure are $t_{b}=t_{t}=0$ and $t_{a}=10 \, \unit{nm}$.
Robin boundary conditions have been enforced 
$\Gamma_{t}$ and $\Gamma_{b}$ 
with $T_{t}=300 \, \unit{K}$ and $T_{b}=900 \, \unit{K}$ and 
$\gamma_{T}=1.17 \cdot 10^{5} \unit{m} \unit{s}^{-1}$, while
homogeneous Neumann conditions are enforced on $\Sigma_{lat}$.
Initial condition for temperature is set constant to $T=300 \, \unit{K}$ 
in all $\Omega_{a}$. For convenience of the reader, 
Tab.~\ref{heat_table} reports the values of the
parameters used during the comparison.
\begin{table}
\begin{center}
\begin{tabular}{|l|l|}
\hline
\textbf{Parameter} & \textbf{value} \\
\hline
$N_{e}$ & $1.0\cdot10^{26} $ \\
$E_{0}$ & $1.158\cdot10^{9}$  \\
$\mu^{el}$ & $3.3\cdot10^{-6} $ \\
$\rho_{m}$ & $3.98\cdot10^{6}$ \\
$c_{m}$ & $880$ \\
$\alpha$ & $10^{-4} $ \\
$\kappa$ & $10^{-1} - 10^{-2}$ \\
\hline 
\end{tabular}
\caption{Parameter values used to test the numerical 
solution of Eq.~\eqref{eq:heat_generalized}.}
\label{heat_table} 
\end{center}
\end{table}

\begin{figure}[h!]
\centering
\includegraphics[width=0.6\textwidth,height=0.5\textwidth]
{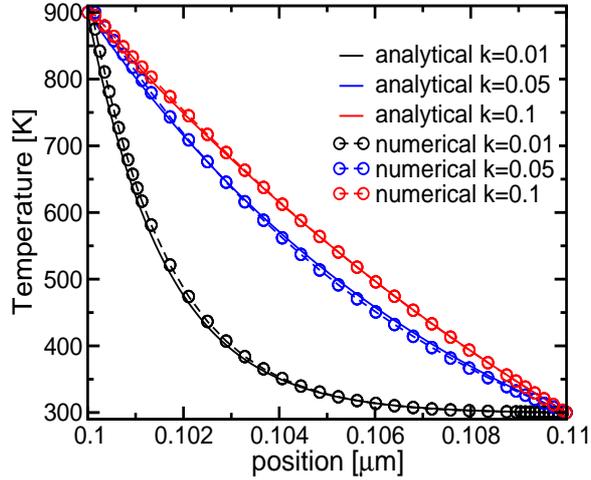}
\caption{Comparison for Eq.~\eqref{eq:heat_generalized} 
between analytical and numerical solution for different values 
of the thermal conductivity. The plot reports for the
3D solution the 1D cut in the center of the $x-y$ plane along 
the $z$ axis.}
\label{fig:thermal_ana}
\end{figure}

Fig.~\ref{fig:thermal_ana} shows a 1D cut along the $z$ 
axis in the center of the $x-y$ plane
of the 3D numerical solutions compared with the exact analytical ones in 
the stationary case (symbols are for numerical and 
lines for analytical values).
To measure the relative weight of thermal diffusion 
with respect to thermal convection it is useful to introduce 
the \emph{local P\`eclet number}
\begin{equation}\label{eq:peclet}
\mathbb{Pe}_{loc} := \Frac{q h \alpha N_{e} \mu_{e}^{el} E_{0}}{2 \kappa}
\end{equation}
$h$ denoting the average mesh size used in the computations, equal
to $10^{-8} \, \unit{m}$. In the three considered cases 
($\kappa=0.01, \, 0.05, \, 0.1 \, \unit{W m^{-1} K^{-1}}$)
the values of $\mathbb{Pe}_{loc}$ are $3$, $0.6$ and $0.3$, respectively,
this indicating that in the first case the thermal flow is dominated by
convection while in the other two cases diffusion is the principal transport
mechanism of heat in the device. It is important to notice that in the
case $\kappa = 0.01 \, \unit{W m^{-1} K^{-1}}$, the use of the exponentially
fitted FEM prevents the onset of spurious oscillations without introducing
any extra amount artificial thermal diffusion.
In all the simulated cases a very good
agreement between numerical and analytical solution is found.

\subsection{Heterogeneous material}\label{sec:heterogeneous}
As discussed in Sect.~\ref{sec:intro}, heterogeneous materials are 
widely used in a new emerging application like PCM devices.
In this section we report the results of 
the numerical simulations in a heterogeneous medium 
for the model of Sect.\ref{sec:model}. 
The test cases (denoted a, b and c) consist in
cubic structures ($t_{b}=t_{t}=0$ and $t_{a}=10 \, \unit{nm}$) where 
the transport region, $\Omega_{a}$, is divided
along the $z$ axis into three zones with thickness of 3, 4 and 
3$\unit{nm}$, respectively. Tab.~\ref{hetero_table_par} 
reports the different parameters used in each of the 
regions separated by a comma.
For sake of clarity, Tab.~\ref{hetero_table_bound} 
shows the boundary conditions applied to the simulation domain
for Eq.~\eqref{eq:Poisson_generalized} to~\eqref{eq:heat_generalized}, 
the symbols D, N and R denoting Dirichlet, Neumann and Robin types,
respectively.
We note that thermal and electrical gradients are 
directed towards $\Gamma_{b}$.
Initial conditions are set constant in $\Omega_{a}$ for
all transport equations to the value of $N_{e}=10^{16} \, \unit{m}^{-3}$ 
for electrons and $T=300 \, \unit{K}$ for temperature.
\begin{table}
\begin{center}
\begin{tabular}{|l|l|l|l|}
 \hline
\textbf{Parameter} & \textbf{a} & \textbf{b} & \textbf{c} \\
 \hline
$\mu^{el}$ & $3\cdot10^{-6},300,3\cdot10^{-10}$& 
$300,300,3\cdot10^{-10}$&$300,3\cdot10^{-10},3\cdot10^{-6}$\\
$\rho_{m}$ & $3.98,3.98,3.98$&$3.98,3.98,3.98$&
$3980,3980,3.98$\\
$c_{m}$ &$880,880,880$&$880,880,880$&$8800,8800,880$\\
$\alpha$ &$10^{-4},10^{-4},10^{-4}$&$10^{-4},10^{-4},10^{-4}$&
$10^{-4},10^{-4},10^{-4}$\\
$\kappa$ &$30,3,300$&$0.3,0.3,300$&$0.03,0.03,300$\\
\hline 
\end{tabular}
\caption{Parameter values used to test the model 
of Sect.~\ref{sec:model} in a heterogeneous structure obtained
dividing $\Omega_{a}$ into three different regions of thickness 
3, 4 and 3$\unit{nm}$ along the $z$ axis.}
\label{hetero_table_par} 
\end{center}
\end{table}

\begin{table}
\begin{center}
\begin{tabular}{|c||l|l|l|}
 \hline
\textbf{Equation} & \textbf{boundary} & \textbf{type} & \textbf{value}\\
 \hline
Eq.~\eqref{eq:Poisson_generalized}& $\Sigma_{b}$ & D & $\varphi=0$\\
Eq.~\eqref{eq:Poisson_generalized}& $\Sigma_{t}$ & D & $\varphi=1$\\
Eq.~\eqref{eq:Poisson_generalized}& $\Sigma_{lat}$ & N & Homogeneous\\
Eq.\eqref{eq:electron_continuity}& $\Gamma_{b}$ & R & 
$v_{eq}=2\cdot10^{2}$; $n_{eq}=10^{19}$\\
Eq.~\eqref{eq:electron_continuity}& $\Gamma_{t}$ & R & 
$v_{eq}=2\cdot10^{2}$; $n_{eq}=10^{13}$\\
Eq.~\eqref{eq:electron_continuity}& $\Sigma_{lat}$ & R & 
$v_{eq}=2\cdot10^{2}$; $n_{eq}=10^{13}$\\
Eq.~\eqref{eq:heat_generalized}& $\Gamma_{b}$ & R & 
$\gamma_{T}=10^{5}$; $T=300$\\
Eq.~\eqref{eq:heat_generalized}& $\Gamma_{t}$ & R & 
$\gamma_{T}=10^{5}$; $T=600$\\
Eq.~\eqref{eq:heat_generalized}& $\Sigma_{lat}$ & R & 
$\gamma_{T}=10^{5}$; $T=300$\\
\hline 
\end{tabular}
\caption{Boundary conditions used to test
the model of Sect.~\ref{sec:model} in the case of a heterogeneous medium.
}
\label{hetero_table_bound} 
\end{center}
\end{table}

\begin{figure}[h!]
\centering
\subfigure[case a]
{\includegraphics[width=0.35\textwidth,height=0.4\textwidth]{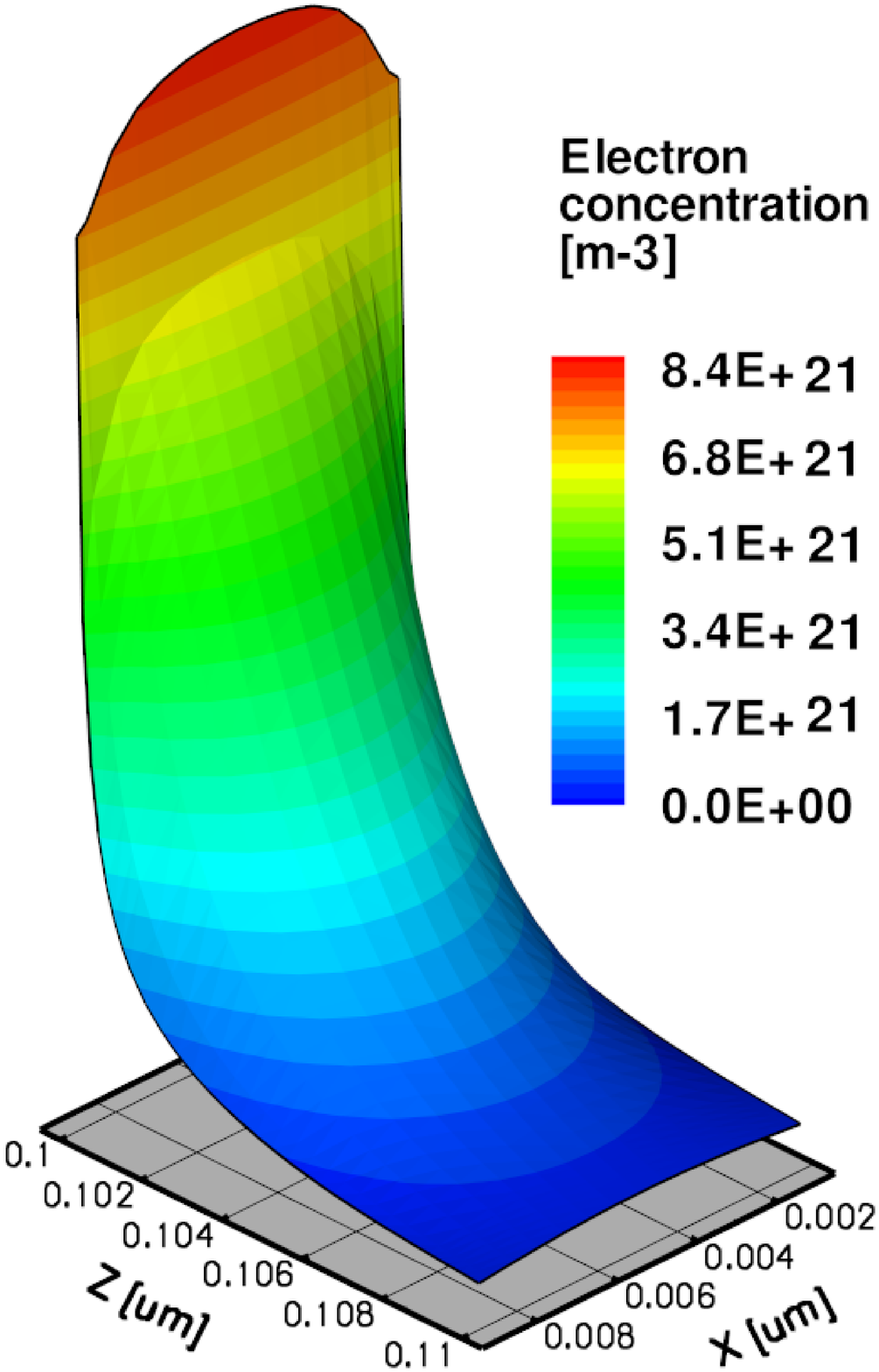}
\label{fig:casea3d_2}}
\subfigure[case b]
{\includegraphics[width=0.35\textwidth,height=0.4\textwidth]{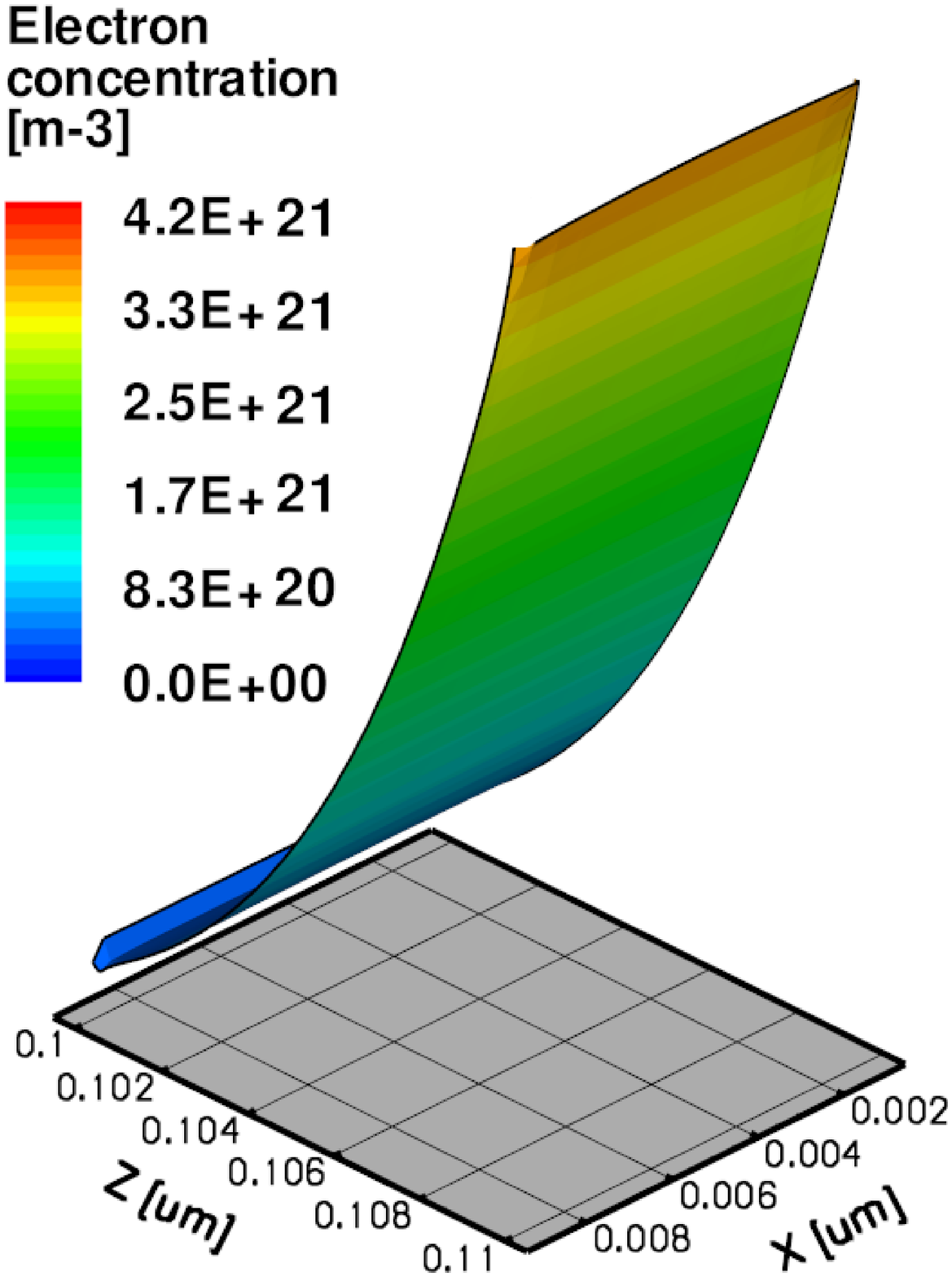}
\label{fig:caseb3d_2}}
\subfigure[case c]
{\includegraphics[width=0.35\textwidth,height=0.4\textwidth]{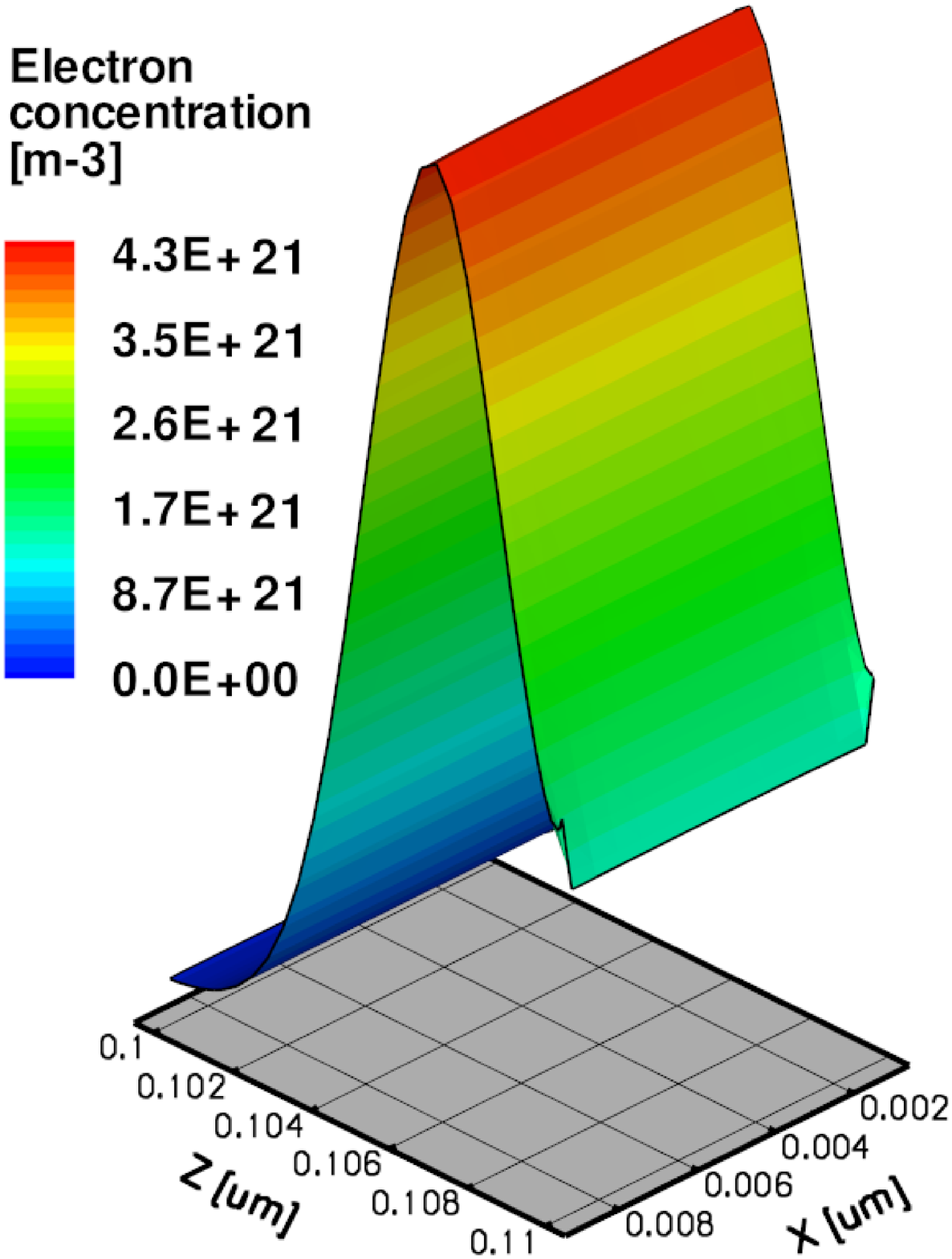}
\label{fig:casec3d_2}}
\caption{Heterogeneous media: electron profiles for cases a), 
b) and c) as reported in Tab.~\ref{hetero_table_par}.}
\label{fig:hetero3d}
\end{figure}

Fig.~\ref{fig:hetero3d} shows the electron concentration obtained 
by the numerical simulations in the three different cases:
because of symmetry reasons we have reported the solutions on the $z$ axis 
and the coordinates of the $x-z$ simulation plane on the $x-y$ plane.

For case a) electrons are forced towards the bottom interface 
place at $z=0.1 \mu m$ because thermal power is high enough 
to force electrons moving against the electric field.
Case b) is exactly the opposite of case a): electrons are moving
along the electric field but against the thermal gradient towards 
$\Gamma_{t}$: the different values of the peaks
for cases a) and b) depend of the different values of 
the electron mobility chosen in the device regions.
More complex to interpret 
are the results showed for case c) in which a charge 
accumulation is found in the center (along the $z$ axis) of the 
active regions:
this is due to the chosen low electron mobility in this region 
and the opposite effects of thermal and electrical gradients.

\begin{figure}[h!]
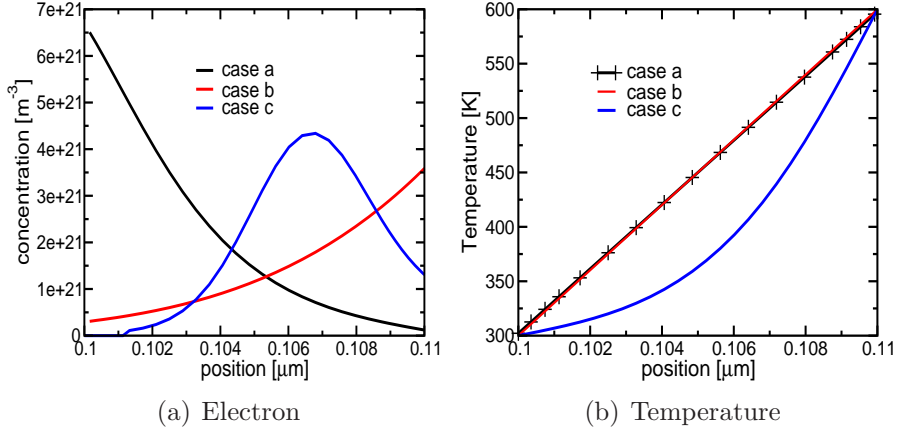

\centering
\subfigure[Electron]
{\includegraphics[width=0.45\textwidth,height=0.4\textwidth]{electron.eps}
\label{fig:heteroel1d}}
\subfigure[Temperature]
{\includegraphics[width=0.45\textwidth,height=0.4\textwidth]{temp.eps}
\label{fig:heterotemp1d}}
\caption{Heterogeneous material: stationary solutions. 
1D cuts in the center of the $x-y$ plane
and along the $z$ axis. Left: Electron concentration. Right: Temperature.}
\label{fig:hetero1d}
\end{figure}

1D cuts along the $z$ axis in the center of the $x-y$ 
simulation plane are shown in Fig.~\ref{fig:hetero1d}. 
In particular, Fig.~\ref{fig:heteroel1d} shows electron 
concentration as in Fig.~\ref{fig:hetero3d} clarifying the accumulation 
of the electrons at the top and bottom interface or at 
the center of the cube.
Fig.~\ref{fig:heterotemp1d} shows 
the temperature profiles: the different 
thermal conductivity chosen for the various region of the devices is resulting
in different thermal velocities 
justifying the difference in the profiles. 

\subsection{Cylindrical shape}\label{sec:cylindrical}
Sect.~\ref{sec:intro} has already pointed how geometries 
can be very complex in the new devices due to 
the miniaturization technological process.
A typical example is the case of a cylindrical shape that can 
simulate the gate all around or three gate devices such those 
employed in the SONOS memories~\cite{sonos2010}, 
or in the FinFET CMOS~\cite{finfet1998} and~\cite{finfet1999} .
The numerical implementation of the model of Sect.~\ref{sec:model} 
has been here applied to study the cylinder of Fig.~\ref{fig:cyl_mesh} 
with $t_{b}=t_{a}=10 \unit{nm}$ and $t_{t}=5 \unit{nm}$.
In Tabs.~\ref{cyl_table_par} and~\ref{cyl_table_bound} 
we have reported the parameters and the boundary conditions
used in the simulations.
We note that thermal and electrostatic gradients
are in the same directions towards the center of the cylinder
in the last two boundary condition for Eq.~\eqref{eq:Poisson_generalized},
while in the first boundary condition they are in the opposite direction.
Moreover the boundary conditions result in an injection
of electrons from the surface $\Sigma_{b}$ inside $\Omega_{a}$.
Initial conditions are set constant in $\Omega_{a}$ 
in the transport equations with $N_{e}=10^{20} \, \unit{m}^{-3}$ 
for electrons and $T=300 \, \unit{K}$ for temperature.
\begin{table}
\begin{center}
\begin{tabular}{|l||l|}
 \hline
\textbf{Parameter} & \textbf{value} \\
 \hline
$\mu^{el}$  & $3\cdot10^{-6}$\\
$\rho_{m}$  & $3.98 \cdot 10^3$\\
$c_{m}$  &    $880$\\
$\alpha$    & $10^{-3}$\\
$\kappa$ & $30$\\
\hline 
\end{tabular}
\caption{Parameter values used to test the model 
of Sect.~\ref{sec:model} in a cylindrical 3D shape.
}
\label{cyl_table_par} 
\end{center}
\end{table}

\begin{table}
\begin{center}
\begin{tabular}{|l||l|l|l|}
 \hline
\textbf{Equation} & \textbf{boundary} & \textbf{type} & \textbf{values}\\
 \hline
Eq.~\eqref{eq:Poisson_generalized}& $\Sigma_{b}$ & D & $\varphi=0$\\
Eq.~\eqref{eq:Poisson_generalized}& $\Sigma_{t}$ & D & 
$\varphi=-0.1,0.4,0.9$\\
Eq.~\eqref{eq:Poisson_generalized}& $\Sigma_{lat}$ & N & Homogeneous\\
Eq.~\eqref{eq:electron_continuity}& $\Gamma_{b}$ & R & 
$v_{eq}=2\cdot10^{2}$; $n_{eq}=10^{25}$\\
Eq.~\eqref{eq:electron_continuity}& $\Gamma_{t}$ & R & 
$v_{eq}=2\cdot10^{2}$; $n_{eq}=10^{19}$\\
Eq.~\eqref{eq:electron_continuity}& $\Sigma_{lat}$ & R & 
$v_{eq}=2\cdot10^{2}$; $n_{eq}=10^{19}$\\
Eq.~\eqref{eq:heat_generalized}& $\Gamma_{b}$ & R & 
$\gamma_{T}=10^{5}$; $T=300$\\
Eq.~\eqref{eq:heat_generalized}& $\Gamma_{t}$ & R & 
$\gamma_{T}=10^{5}$; $T=600$\\
Eq.~\eqref{eq:heat_generalized}& $\Sigma_{lat}$ & R & 
$\gamma_{T}=10^{5}$; $T=300$\\
\hline 
\end{tabular}
\caption{Boundary conditions used to test
the model of Sect.~\ref{sec:model} in a cylindrical 3D shape.
}
\label{cyl_table_bound} 
\end{center}
\end{table}

\begin{figure}[h!]
\centering
\subfigure[$V_{top} =-0.1V$]
{\includegraphics[width=0.3\textwidth,height=0.3\textwidth]{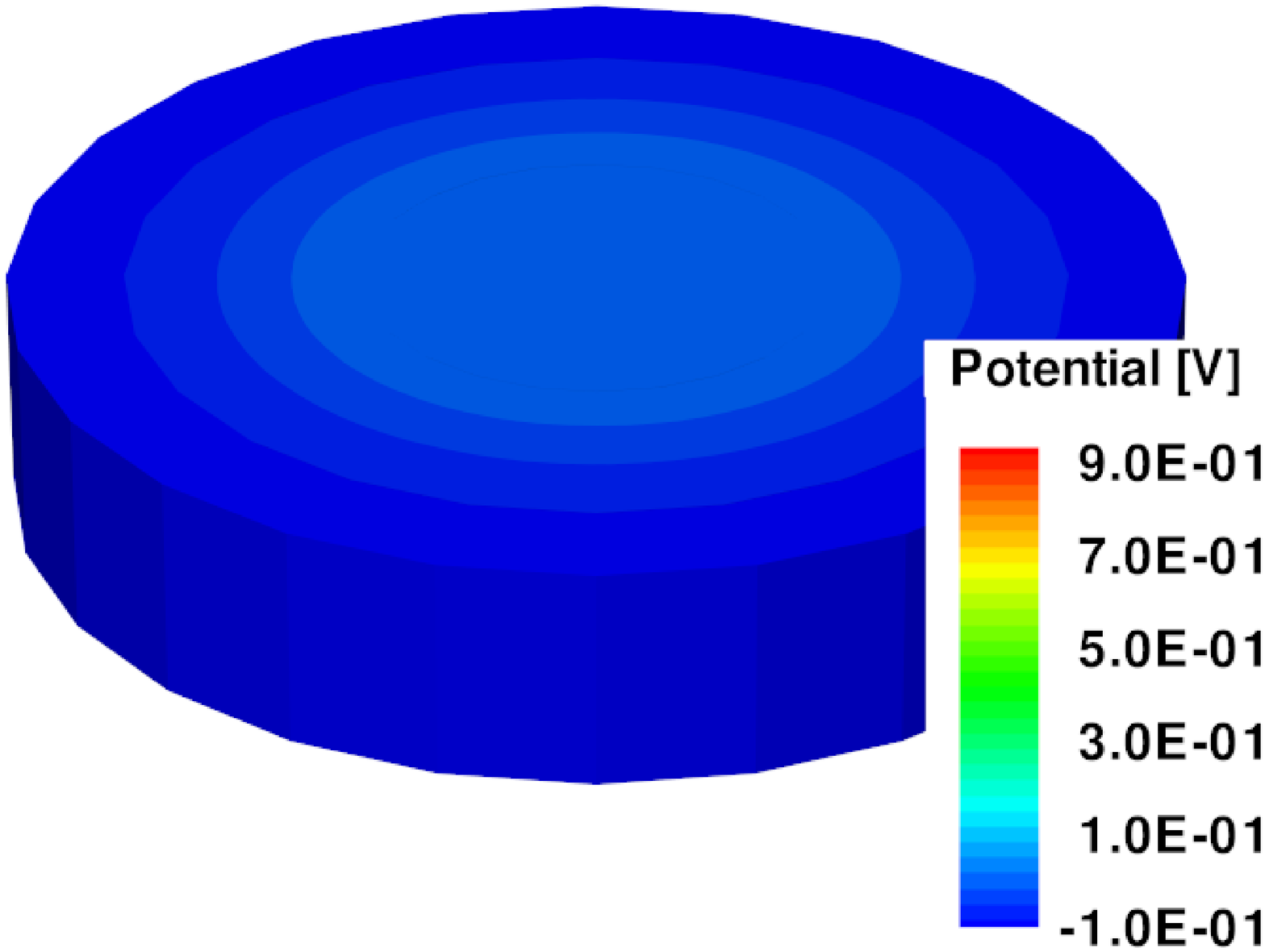}
\label{fig:casea3d}}
\subfigure[$V_{top} =0.4V$]
{\includegraphics[width=0.3\textwidth,height=0.3\textwidth]{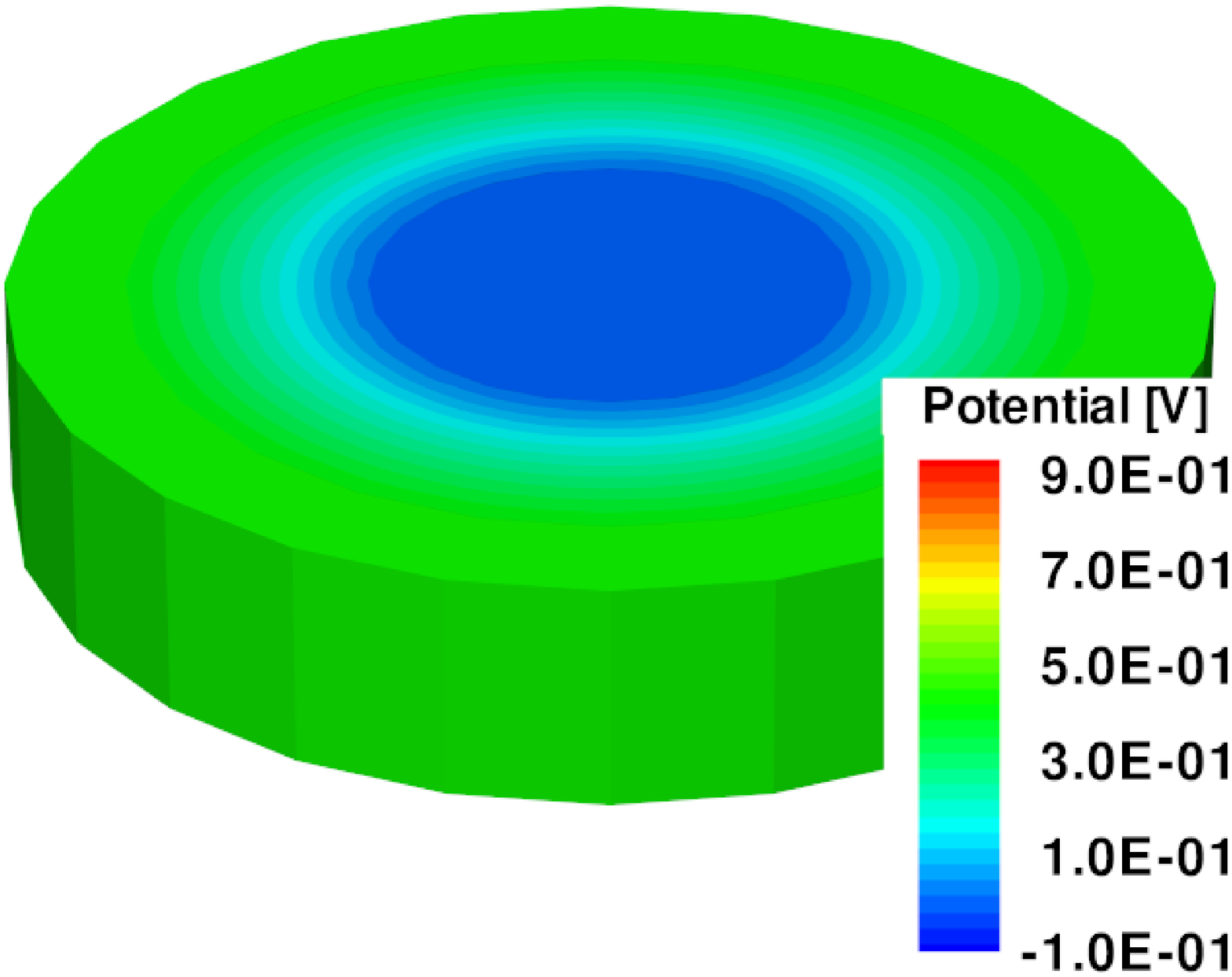}
\label{fig:caseb3d}}
\subfigure[$V_{top} =0.9V$]
{\includegraphics[width=0.3\textwidth,height=0.3\textwidth]{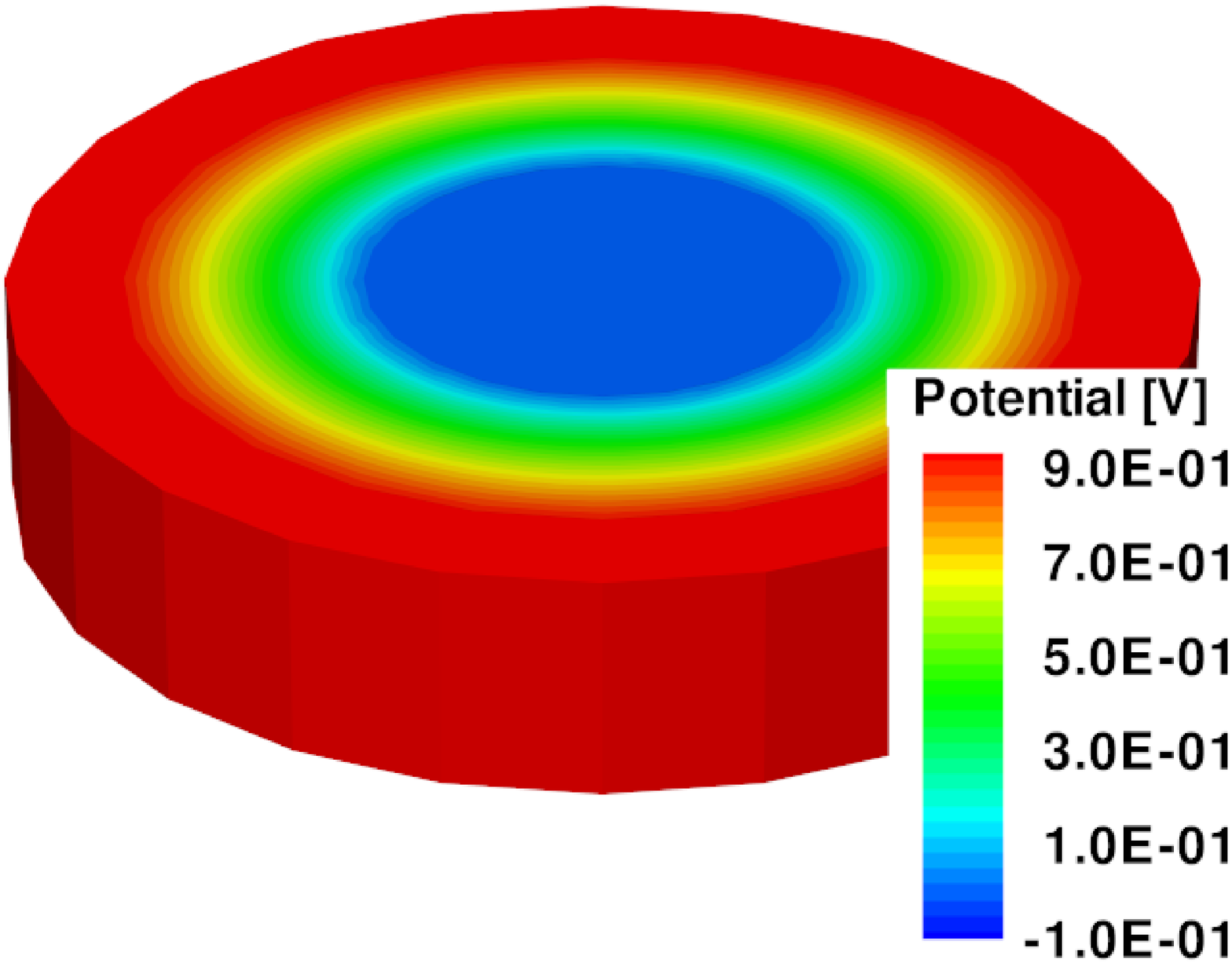}
\label{fig:casec3d}}
\caption{Cylindrical shape: Solution of 
Eq.~\eqref{eq:Poisson_generalized}
for different Dirichlet boundary conditions.}
\label{fig:cyl3d}
\end{figure}

Fig.~\ref{fig:cyl3d} shows the numerical solution of 
Eq.~\eqref{eq:Poisson_generalized}
for the three different applied bias on $\Sigma_{t}$: as expected, 
the potential is a continous function overall the device 
and the gradient direction is swapping between 
the first and the last two values.

\begin{figure}[h!]
\centering
\includegraphics[width=0.7\textwidth,height=0.5\textwidth]{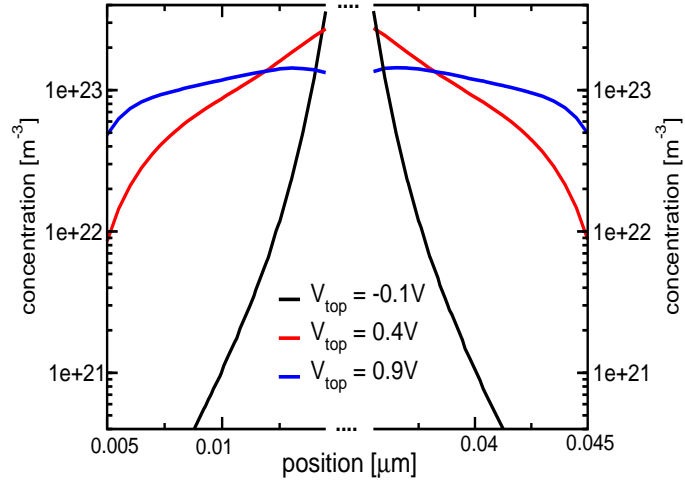}
\caption{1D cut at $z=0.005 \mu m$ and $y=0.025 \mu m$ of the 
3D numerical solution of Eq.~\eqref{eq:electron_continuity} 
for different applied potentials.}
\label{fig:cyl_1d}
\end{figure}

Fig.~\ref{fig:cyl_1d} shows 1D cuts at $z=0.005 \mu m$ and $y=0.025 
\mu m$ of electron concentration (cf. Eq.~\eqref{eq:electron_continuity})
for three different applied bias on $\Sigma_{t}$.  
Even if electrons are injected in $\Omega_{a}$ 
from $\Sigma_{b}$, in the case with $V=-0.1 V$ on $\Sigma_{b}$ 
electric and thermal gradients are pushing electrons 
back towards $\Sigma_{b}$. This results in very low diffused profiles.
In the latter two cases, on the contrary, the electric field is high 
enough to dominate over the thermal gradient and hence 
electrons can diffuse towards $\Sigma_{t}$.

\section{Conclusions and Future Perspectives}\label{sec:conclusions}

In the present article we have provided a unified 
mathematical framework capable of describing the complex 
and interplaying electro-thermo-chemical
processes that occur in modern new emerging technologies in semiconductor
device industry. 

The general conservation law format of the model building block
equations allowed us to successfully adapt to the presently investigated
application: 1) the functional iteration tools usually employed in 
standard electronic transport device simulation programs, 
and 2) the Finite Element Exponentially 
Fitted discretization technique that, in conjunction with a suitable 
tetrahedral geometrical partition of the computational domain, is
characterized by enjoying a discrete maximum principle for chemical
number densities and temperatures.

Model and computational algorithms have been thoroughly validated 
by the numerical study of several realistic device geometries for 
which, in some simple albeit significant cases, exact analytical 
solutions were available. Results have always been characterized by 
a very good mathematical accuracy and close agreement with physically 
expected solution behaviour, clearly demonstrating the potentiality of 
model and numerical tools in providing close insights and fine
prediction for outperforming devices of the next node generation.

Future steps in our research programme in this new emerging area include:
\begin{enumerate}
\item further validation of the proposed computational model
through calibration against measured data;
\item inclusion of mechanical stress analysis in the model;
\item analysis of the existence of a fixed point and convergence 
of the functional iteration~\eqref{eq:gummel_thec};
\item analysis of well-posedness of each differential
subproblem in the iterative map~\eqref{eq:gummel_thec};
\item proof of local/global estimates in time of the solution
of the PDE system~\eqref{eq:model_thec}-~\eqref{eq:fluxes_ohmic_diel_thec} 
supplemented by the initial conditions~\eqref{eq:ICs} 
and boundary conditions~\eqref{eq:BCs}.
\end{enumerate}

\section*{Acknowledgements}\label{sec:ack}
The authors gratefully acknowledge 
Giovanni Novielli and Silvia Sorbello 
(MSc degree students in Mathematical Engineering
at Politecnico di Milano) for their contribution in the development
of the simulation program and of the numerical results.

\bibliographystyle{amsplain}
\bibliography{biblio_ram}

\appendix

\section{List of symbols}\label{sec:list}

Below, we provide a summary of all the variables,
physical constants and parameters that have been introduced
throughout the article, specifying for each symbol the meaning
and units.

\begin{table}[h!]
\begin{center}
\begin{tabular}{llll}
{\bf Symbol}         & {\bf Meaning} & {\bf Units} \\[3mm]
$\vect{x}$ & position vector & $\meter$ \\
$t$        & time variable & $\second$ \\
$M$  & total number of chemicals & \\
$z_i$ & ionic valence & \\
$N_i$ & number density & $\meter^{-3}$ \\
$T$ & temperature & $\unit{K}$ \\
$\varphi$ & electric potential & $\volt$ \\
$\vect{E}$ & electric field & $\volt \meter^{-1}$ \\
$R_i$ & net production rate & $\meter^{-3} \second^{-1}$ \\
$\vect{j}_i$ & current density & $\unit{A} \meter^{-2}$ \\
$\vect{j}_T$ & energy flux density & $\unit{W} \meter^{-2}$ \\
$\rho$ & mass density & $\unit{Kg} \, \meter^{-3}$ \\
$c$ & specific heat & $\meter^2 \second^{-2} \unit{K}^{-1}$ \\
$\mathcal{Q}_T$ & net heat production rate & $\unit{W} \meter^{-3}$ \\
$\varepsilon$ & dielectric permittivity & 
$\unit{F} \, \meter^{-1}$ \\
$\mathcal{D}$ & net doping & $\meter^{-3}$ \\
$\vect{j}_i^{ec}$ & electrochemical current flux & $\unit{A} \meter^{-2}$ \\
$\vect{j}_i^{th}$ & thermal current flux & $\unit{A} \meter^{-2}$ \\
$\vect{j}_T^{ec}$ & electrochemical heat flux & $\unit{W} \meter^{-2}$ \\
$\vect{j}_T^{th}$ & thermal heat flux & $\unit{W} \meter^{-2}$ \\
$\mu_i^{el}$ & electrical mobility & $\meter^2 \volt^{-1} 
\second^{-1}$ \\
$\sigma_i$ & electrical conductivity & $\unit{S} \meter^{-1}$ \\
$\varphi_i^{c}$ & chemical potential & $\volt$ \\
$\varphi_i^{ec}$ & electrochemical potential & $\volt$ \\
$\mu_i^c$ & chemical energy & $\unit{J} \unit{mol}^{-1}$ \\
$N_{ref}$ & reference concentration & $\meter^{-3}$ \\
$\varphi^{th}$ & thermal potential & $\volt$ \\
$\alpha$ & thermopower coefficient & $\volt \unit{K}^{-1}$ \\
$\vect{E}_i^{thec}$ & thermo-electrochemical field & $\volt \meter^{-1}$ \\
$\psi_i$ & thermo-electrochemical potential & $\volt$ \\
$D_i$ & diffusion coefficient & $\meter^2 \second^{-1}$ \\
$\vect{E}_i^{el}$ & generalized electric field & $\volt \meter^{-1}$ \\
$\kappa$ & thermal conductivity & $\unit{W} \meter^{-1} \unit{K}^{-1}$ \\
$\vect{j}$ & total thermo-electrochemical flux & $\unit{A} \meter^{-2}$ \\
$\psi$ & total thermo-electrochemical potential & $\volt$ \\
\end{tabular}
\label{tab:model_params} 
\end{center}
\end{table}

\begin{table}[h!]
\begin{center}
\begin{tabular}{llll}
{\bf Symbol}         & {\bf Meaning} & {\bf Units} & 
{\bf Numerical Value} \\[3mm]
$q$ & electron charge & $\unit{C}$ & $1.602 \cdot 10^{-19}$ \\
$K_B$ & Boltzmann constant & $\unit{J} \unit{K}^{-1}$ & 
$1.38 \cdot 10^{-23}$ \\
$R$   & ideal gas constant & 
$\unit{J} \unit{K}^{-1} \unit{mol}^{-1}$ & $8.314$ \\
$F$ & Faraday constant & 
$\unit{C} \unit{mol}^{-1}$ & $9.648 \cdot 10^4$ \\
\end{tabular}
\label{tab:constants} 
\end{center}
\end{table}

\end{document}